\input amstex
\documentstyle{amsppt} 
\TagsOnRight
\magnification=1200
\hcorrection{.25in}
\advance\vsize-.75in
\overfullrule0pt
\NoBlackBoxes

\def\block #1{\vbox{\hsize 2.5 true cm\noindent \bf#1:}}

\def\memo TO:#1FROM:#2SUBJECT:#3DATE:#4\par{\centerline{\sl MEMO!}
      \bigskip \hrule height1pt \medskip
      \vbox{\parindent=75pt\parskip=1pt
      \item{\block{TO}}#1
      \item{\block{FROM}}#2
      \item{\block{SUBJECT}}#3
      \item{\block{DATE}}#4 }\medskip \hrule height1pt \bigskip}

%********** SCRIPTS ************************************************************
\newcount\spk
\def\beginscript {\bgroup \parindent=0pt \spk=1 \sl \rightskip.4in
      \def\par {\ifnum\spk=1 \endgraf \it \spk=2 \leftskip.4in \rightskip0in
               \else \endgraf \sl \spk=1 \leftskip0in \rightskip.4in \fi}}

\def\endscript {\egroup}

%********** FERMAT *************************************************************
\newcount\var \newcount\pw \newcount\tmp \newcount\cnt
\def\pow#1#2#3{\var=#1 \pw=#2 \tmp=\var \cnt=1
      \loop \multiply\var by\tmp \advance\cnt by1 \ifnum\cnt<\pw \repeat
      \global#3=\var}

\newcount\xf \newcount\xnf \newcount\yf \newcount\ynf \newcount\zf \newcount\znf
\def\n {\number}

\def\fermat#1#2#3{$\global\xf=#1 \global\yf=#2 \global\pw=#3
      \pow{\xf}{\pw}{\xnf} \pow{\yf}{\pw}{\ynf}
      \global\tmp=\ynf \global\advance\tmp by\xnf
      {\n\xf}^{\n\pw}+{\n\yf}^{\n\pw}={\n\tmp}$.\hfil\break
      \ifnum\xf>\yf \zf=\xf \else \zf=\yf \fi
      \loop {\pow{\zf}{\pw}{\znf}} \ifnum\znf<\tmp \advance\zf by1 \repeat
      \ifnum\znf=\tmp  The sum seems to be exactly ${\n\zf}^{\n\pw}$.
        \ifnum\pw=2 {\it Yawn!\/} Tell me something I don't know, will you?
        \else Incredible!  But, perhaps you'd better check my work.\fi
      \else \advance\zf by-1
      This lies between ${\n\zf}^{\n\pw}$\pow{\zf}{\pw}{\znf}($={\n\znf}$)
      \advance\zf by1
      and ${\n\zf}^{\n\pw}$\pow{\zf}{\pw}{\znf}($={\n\znf}$).\fi}%

%********** FRAMING ************************************************************
\def\frame #1#2#3#4{\vbox{\hrule height #1pt%    TOP RULE
 \hbox{\vrule width #1pt\kern #2pt%              RULE AND SPACE ON LEFT
 \vbox{\kern #2pt%                               SPACE AT TOP
 \vbox{\hsize #3\noindent #4}%                   MATERIAL THAT WILL BE BOXED
 \kern #2pt}%                                    SPACE AT BOTTOM
 \kern #2pt\vrule width #1pt}%                   SPACE AND RULE ON RIGHT
 \hrule height0pt depth #1pt}}%,                 BOTTOM RULE

\def\fitframe #1#2#3{\vbox{\hrule height#1pt%    TOP RULE
 \hbox{\vrule width#1pt\kern #2pt%               RULE AND SPACE ON LEFT
 \vbox{\kern #2pt\hbox{#3}\kern #2pt}%           TOP, MATERIAL, BOTTOM
 \kern #2pt\vrule width#1pt}%                    SPACE AND RULE ON RIGHT
 \hrule height0pt depth#1pt}}%                   BOTTOM RULE

\def\shframe #1#2#3#4{\vbox{\hrule height 0pt%   NO TOP SHADOW
 \hbox{\vrule width #1pt\kern 0pt%               LEFT SHADOW
 \vbox{\kern-#1pt\frame{.3}{#2}{#3}{#4}%         SHADOW STARTS #1 PT FROM TOP
 \kern-.3pt}%                                    MOVE UP RULE THICKNESS AT BOT.
 \kern-#1pt\vrule width 0pt}%                    STOPS #1 PT FROM RT; NO RT SHAD
 \hrule height #1pt}}%                           BOTTOM SHADOW

\def\s #1{\frame{.3}{2}{8pt}{\centerline{#1\vphantom{(}}}\ }

%********** RULED TRAPEZOID ****************************************************
\newcount\tw        %top width in pt
\newcount\bw        %bottom width in pt
\newcount\h         %height in pt
\newcount\bs        %bottom shift in pt
\newcount\th        %line thickness in 1/64th of a pt
\newcount\gp        %line gap in 1/64th of a point
\newcount\rs        %running shift
\newcount\rw        %running width
\newcount\rh        %running height
\newcount\tmp       %for temporarily storing variables

\def\trap #1#2#3#4#5#6{\vbox{\offinterlineskip
      \tw=#1 \bw=#2 \h=#3 \bs=#4 \th=#5 \gp=#6 \rh=0
      \multiply\tw by 65536 \multiply\bw by 65536 \multiply\bs by 65536
      \multiply\th by 1024 \multiply\gp by 1024
      \loop
      \tmp=\bs \multiply\tmp by\rh \divide\tmp by\h
      \rs=\tmp                                        %running shift calculated
      \tmp=\bw \advance\tmp by-\tw \multiply\tmp by\rh
      \divide\tmp by\h \advance\tmp by\tw \rw=\tmp    %running width calculated
      \hbox{\kern\rs sp\vrule height0sp depth\th sp width\rw sp}%RULE DRAWN HERE
      \vskip\gp sp                                              %GAP LEFT HERE
      \ifnum\rh<\h  \tmp=\rh \multiply\tmp by 65536
       \advance\tmp by\th \advance\tmp by\gp \divide\tmp by65536 \rh=\tmp
      \repeat}}%

%

%********** MULTIPLE COLUMN ****************************************************
\newcount\colnumber \newbox\col \newdimen\tmpdim \newdimen\size
\newdimen\coljump \coljump=.2 true in                  %GAP BETWEEN COLUMNS

\def\nstrut {\vrule height\topskip depth0pt width0pt}  %TO PROP THINGS UP
\def\divider{\hskip\coljump}
\def\dividerule{\dimen0=.4pt \dimen1=\coljump \advance\dimen1 by-\dimen0
      \divide\dimen1 by2
      \def\divider{\hskip\dimen1 \vrule width\dimen0 \hskip\dimen1}}

\def\beginart #1/#2{\vbox\bgroup#1 \colnumber=#2  \parskip=0pt
      \advance\colnumber by-1 \tmpdim=\coljump \multiply\tmpdim by\colnumber
      \size=\hsize  \advance\size by-\tmpdim
      \advance\colnumber by1 \divide\size by\colnumber
      \vbadness=10000 \hbadness=2000 \tolerance=2000
      \setbox\col=\vbox\bgroup\hsize\size \noindent\nstrut}

\def\endart {\global\size=\baselineskip \vfil \egroup
      \multiply\size by\colnumber  \advance\size by-\baselineskip
      \tmpdim=\ht\col  \advance\tmpdim by\size  \divide\tmpdim by\colnumber
      \hbox{\splittopskip=\topskip \doittoit}\egroup}

\def\doittoit{\ifnum\colnumber>0 \vsplit\col to \tmpdim
      \global\advance\colnumber by-1
      \ifnum\colnumber>0 \divider \fi \doittoit \fi}

%********** FONT AND CODE TABLE ************************************************
\def\chartable #1{\smallbreak\vbox{\noindent%
     \underbar{Characters in the {\sl#1\/} font, with decimal codes:}
     \raggedright \hbadness5000 \tolerance10000 \medskip
     \font\ft=#1 \ft \dimen0=14pt \baselineskip=\dimen0
     \ifdim\dimen0<3.25ex \baselineskip=3.25ex \fi \count255=0 \dimen0=10pt
     \loop \setbox0=\hbox{\char\count255}
     \ifdim\wd0>\dimen0 \dimen0=\wd0 \fi \advance\count255 by1
     \ifnum\count255<128 \repeat \count255=0
     \advance\dimen0 by 25pt \noindent \loop
     \hbox to\dimen0{\hbox to23pt{\hfil\rm\the\count255:\ }\char\count255\hfil}
     \advance\count255 by1 \ifnum\count255<128 \quad\repeat
     \smallbreak} \font\tenrm=cmr10 }

%********** THE `IGNORE' COMMAND ***********************************************
\def\ignore {\count255=0 \begingroup
      \loop \catcode\count255=14  % Make everything a comment character.
         \advance\count255 by1 \ifnum\count255<127
      \repeat \catcode`\!=0 }     % Makes ! an escape character.
{\catcode`\!=0 !gdef!E{!endgroup}}% Defines the `stop ignoring' command.

%********** STACKSYMBOLS *******************************************************
\def\stacksymbols #1#2#3#4{\def\theguybelow{#2}
    \def\verticalposition{\lower#3pt}
    \def\spacingwithinsymbol{\baselineskip0pt\lineskip#4pt}
    \mathrel{\mathpalette\intermediary#1}}
\def\intermediary#1#2{\verticalposition\vbox{\spacingwithinsymbol
      \everycr={}\tabskip0pt
      \halign{$\mathsurround0pt#1\hfil##\hfil$\crcr#2\crcr
               \theguybelow\crcr}}}

%********** DATING *************************************************************
\def\monthname {\ifcase\month\or January\or February\or March\or
April\or May\or June\or July\or August\or September\or October\or
November\or December\fi}

\newcount\mins  \newcount\hours  \hours=\time \mins=\time
\def\now{\divide\hours by60 \multiply\hours by60 \advance\mins by-\hours
     \divide\hours by60         % NOTE: \divide only gives integer answers.
     \ifnum\hours>12 \advance\hours by-12
       \number\hours:\ifnum\mins<10 0\fi\number\mins\ P.M.\else
       \number\hours:\ifnum\mins<10 0\fi\number\mins\ A.M.\fi}

%********** TABLE COMMANDS *****************************************************
\newdimen\tempdim                 % For temporary storage.
\newdimen\othick   \othick=.4pt   % To set the outer rule thickness.
\newdimen\ithick   \ithick=.4pt   % To set the inner rule thickness.
\newdimen\spacing  \spacing=9pt   % To set the interline spacing.
\newdimen\abovehr  \abovehr=6pt   % Space above horizontal rules.
\newdimen\belowhr  \belowhr=8pt   % Space below horizontal rules.
\newdimen\nexttovr \nexttovr=1pt  % Space next to vertical
                                  % rules. (Diane change from 8 to 6pt)

\def\r{\hfil&\omit\vrsp\vrule width\othick\cr&}   % To start a new line.
\def\rr{\hfil\down{\abovehr}&\omit\vrsp\vrule width\othick\cr
     \noalign{\hrule height\ithick}\up{\belowhr}&}% To draw an \hrule.
\def\up#1{\tempdim=#1\advance\tempdim by1ex
     \vrule height\tempdim width0pt depth0pt}%   For space above a line.
\def\down#1{\vrule height0pt depth#1 width0pt}%  For space below a line.
\def\large#1#2{\setbox0=\vtop{\hsize#1 \lineskiplimit=0pt \lineskip=1pt
     \baselineskip\spacing \advance\baselineskip by 3pt \noindent
     #2}\tempdim=\dp0\advance\tempdim by\abovehr\box0\down{\tempdim}}
% `\large' allows you to make multi-line table entries.
  % To leave a space, the width of a digit.
 % For a gap.
\def\vrsp{\hskip\nexttovr\relax}
\def\toprule#1{\def\startrule{\hrule height#1\relax}} % Set a top rule.
\toprule{\othick}                      % Picking the `\toprule' default.
\def\nstrut{\vrule height\spacing depth3.5pt width0pt}
                 % To print an exclamation mark.
\def\preamble#1{\def\startup{#1}}      % For `customized' preambles.
\preamble{&##}                         % Choosing the default preamble.
{\catcode`\!=\active
 \gdef!{\hfil\vrule width0pt\vrsp\vrule width\ithick\relax\vrsp&}}
% Setting up `!' as the entry separator.

\def\table #1{\vbox\bgroup \setbox0=\hbox{#1}
     \vbox\bgroup\offinterlineskip  \catcode`\!=\active
     \halign\bgroup##\vrule width\othick\vrsp&\span\startup\nstrut\cr
     \noalign{\medskip}
     \noalign{\startrule}\up{\belowhr}&}

\def\caption #1{\down{\abovehr}&\omit\vrsp\vrule width\othick\cr
     \noalign{\hrule height\othick}\egroup\egroup \setbox1=\lastbox
     \tempdim=\wd1 \hbox to\tempdim{\hfil \box0 \hfil} \box1 \smallskip
     \hbox to\tempdim{\advance\tempdim by-20pt\hfil\vbox{\hsize\tempdim
     \noindent #1}\hfil}\egroup}

%********** VERBATIM REPRODUCTION **********************************************
\let\cc=\catcode
{\cc`\^^M=\active %
\gdef\losenolines{\cc`\^^M=\active \def^^M{\leavevmode\endgraf}}}
\def\literal {\begingroup \cc`\\=12 \cc`\{=12 \cc`\}=12 \cc`\$=12 \cc`\&=12
 \cc`\#=12 \cc`\%=12 \cc`\~=12 \cc`\_=12 \cc`\^=12 \cc`\*=12 \cc`\@=0
 \cc`\`=\active \losenolines \obeyspaces \tt}%
{\obeyspaces\gdef {\hglue.5em\relax}}

{\cc`\`=\active \gdef`{\relax\lq}}

\def\vquotingon{\cc`\"=\active}
\def\vquotingoff{\cc`\"=12}
\vquotingon
\def"{\literal\leavevmode\hbox\bgroup\com}
%`\leavevmode' starts a new paragraph, if needed.

{\cc`\@=0 \cc`\\=12 @cc`@^^M=@active %
 @gdef@com#1"{#1@egroup@endgroup} %
 @gdef@thatisit^^M#1\endliteral{#1@endgroup@smallskip}}
\vquotingoff

%********** GENERAL ************************************************************
\def\pattern #1#2{\count0=0
      \loop #1\advance\count0 by 1 \ifnum\count0<#2 \repeat}

%*******************************************************************************

\def \n {\noindent}
\def \s {\smallskip}
\def \m {\medskip}
\def \b {\bigskip}
\def \a {\alpha}
  
\def \be {\beta}

\def \e {\epsilon} 
\def \End {{\text {End}}}

\def\eqa{\buildrel\text {1} \over =} 
\def\eqb{\buildrel\text {2} \over =}
\def\eqc{\buildrel\text {3} \over =}

\def \ga {\gamma}  
\def \g {{\frak g}} 
\def \h {{\frak h}} 
\def \Hom {{\text {Hom}}}
\def \l {\lambda}

\def\map1{\buildrel{p \text{is a}\uqd-\text{hom}} \over =} 
\def \o {\overline}
\def \ot {\otimes}

\def\fimap {\buildrel{ (\,|\,) \ot \text{id}} \over \rightarrow}
\def\jimap {\buildrel{ \jmath \ot \text{id}} \over \rightarrow} 
\def\jpmap {\buildrel{ \jmath \circ p \circ (\jmath \ot \jmath)} \over \rightarrow}
\def\pmap {\buildrel  p \over \rightarrow}
\def\famap {\buildrel  F_{\a_3}+F_{\a_4} \over \rightarrow}
\def\fbmap {\buildrel  F_{\a_2} \over \rightarrow}
\def\fcmap {\buildrel  F_{\a_1} \over \nearrow} 
\def\fdmap {\buildrel  F_{\a_3}+F_{\a_4} \over \searrow}
\def\lmap {\buildrel{ \text{id}\ot\Big( \jmath 
\circ p^* \circ (\jmath \ot \jmath)\Big)} \over \longleftarrow}
\def \ppmap {\buildrel  \text{id} \ot p^* \over \rightarrow} 

\def \u {\underline}
\def \uqd {U_q(\text{D}_4)}

\def \F {{\text {\bf F}}}

\def \K {{\text {\bf K}}} 
\def \Z {{\text {\bf Z}}}

\topmatter
\title  A Quantum Octonion Algebra 
\endtitle
\author Georgia Benkart\\
Jos\'e M. P\'erez-Izquierdo
\endauthor
\address
\hskip-\parindent
Department of Mathematics, University of Wisconsin, Madison,
Wisconsin 53706, USA 
\endaddress
\thanks The first author gratefully acknowledges support from National Science
Foundation Grant \#{}DMS--9622447.
\newline \indent
The second author is grateful for support from the Programa
	de Formaci\'on del Personal Investigador en el Extranjero and
	from Pb 94-1311-C03-03, DGICYT.
\newline \indent
Both authors acknowledge with gratitude the
support and hospitality of the Mathematical
Sciences Research Institute, Berkeley.  
Research at MSRI is supported in part by NSF grant DMS-9701755.
\endthanks

\subjclass Primary 17A75, 17B37, 81R50\endsubjclass

\email
benkart\@math.wisc.edu, perez\@math.wisc.edu
\endemail

%\leftheadtext{GEORGIA BENKART, JOS\'E M. P\'EREZ-IZQUIERDO}
%\rightheadtext {A QUANTUM OCTONION ALGEBRA} 

\endtopmatter 

\centerline {\it To the memory of Alberto Izquierdo}

\abstract
Using the natural irreducible 8-dimensional representation 
and the two spin representations of the quantum group $U_q$(D$_4$) of
D$_4$,  we construct a quantum  analogue  of the split octonions and
study its properties. 
 We  prove that 
the quantum octonion algebra  satisfies  the q-Principle of Local Triality
and has a nondegenerate bilinear
form which satisfies a q-version of the composition
property.   By its construction, the quantum octonion algebra is
a nonassociative algebra with a Yang-Baxter operator action coming from 
the R-matrix of $U_q$(D$_4$). 
The product in the quantum octonions  is a $U_q$(D$_4$)-module
homomorphism.  Using that, we prove identities
for the quantum octonions, and as a consequence, obtain at  $q = 1$
new ``representation
theory'' proofs for very well-known identities satisfied by the octonions.
In the process of constructing the quantum octonions
we introduce an algebra which is a q-analogue of the 
8-dimensional para-Hurwitz algebra.
\endabstract

\head {Introduction} \endhead 
Using the representation theory of the quantum group
$U_q$(D$_4$) of D$_4$  we construct a quantum analogue $\Bbb O_q$ of
the split octonions and study its properties.  
\m 
A unital algebra over a field with a nondegenerate
bilinear form $(\,|\,)$ of maximal Witt index which admits composition, 

$$(x\cdot y | x \cdot y) = (x|x)(y|y),$$

\n must be the field, two copies of the field, the split quaternions,
or the split octonions.   There is a natural $q$-version
of the composition property that the algebra $\Bbb O_q$ of quantum octonions
is shown to satisfy (see Prop. 4.12 below).   We also prove that 
the quantum octonion algebra $\Bbb O_q$ 
satisfies the ``q-Principle of Local Triality'' (Prop. 3.12).  Inside the quantum
octonions are two nonisomorphic 4-dimensional subalgebras, which
are $q$-deformations of the split quaternions.  One of them is unital,
and both of them give $gl_2$ when considered as algebras under the
commutator product $[x,y] = x\cdot y - y \cdot x$.  
\m   
By its construction, $\Bbb O_q$ is
a nonassociative algebra with a Yang-Baxter operator action coming from 
the R-matrix of $\uqd$.  
Associative algebras with a Yang-Baxter operator (or $r$-algebras) arise in Manin's
work on noncommutative geometry and include such important
examples as Weyl and Clifford algebras, quantum groups,
and certain universal enveloping algebras (see for example [B]).
\m
In the process of constructing $\Bbb O_q$ we define an algebra $\Bbb P_q$,
which is a $q$ analogue of the 8-dimensional para-Hurwitz algebra.
The quantum para-Hurwitz algebra  $\Bbb P_q$ is shown to
satisfy certain identities which become familiar properties of
the para-Hurwitz algebra at $q = 1$.  The para-Hurwitz algebra is
a (non-unital) 8-dimensional algebra with a non-degenerate 
associative bilinear form 
admitting composition.  The algebra
$\Bbb P_q$ exhibits related properties (see Props. 5.1, 5.3). 
\m

While this paper was in preparation, we received a preprint of
[Br], which constructs a quantized octonion algebra using the
representation theory of $U_q(sl_2)$.  Although both
Bremner's quantized octonions and our quantum octonions reduce
to the octonions at $q = 1$, they are different algebras.
The quantized octonion algebra in [Br] is constructed
from defining a multiplication on the sum of 
irreducible $U_q(sl_2)$-modules of dimensions 1 and 7.
As a result, it carries a $U_q(sl_2)$-module structure and has
a unit element.  
The quantum octonion algebra constructed in this
paper using $\uqd$  has a unit element only for the
special values $q = 1,-1$. 
\m

An advantage to the  $\uqd$ approach is that allows 
properties such as the ones mentioned above to be derived from
its representation theory.  Using
the fact that the product in $\Bbb O_q$  is a $\uqd$-module
homomorphism, we prove identities
for $\Bbb O_q$ (see Section 4)  and as a consequence obtain at  $q = 1$ new ``representation
theory'' proofs for very well-known identities satisfied by the octonions
that had been established previously  by other methods.   
The $\uqd$ approach also affords connections with fixed points of
graph automorphisms - although as we show in Section 7, the
fixed point subalgebras of $\uqd$ are not the
quantum groups $U_q$(G$_2$)
and $U_q($B$_3)$, because those quantum
groups do not have Hopf algebra embeddings into $\uqd$. 
In the final section of this paper we explore connections
with quantum Clifford algebras.  In particular, we obtain a $\uqd$-isomorphism
between the quantum Clifford algebra $C_q(8)$
and the endomorphism algebra $\End(\Bbb O_q \oplus \Bbb O_q)$ and discuss its
relation to the work of Ding and Frenkel [DF] (compare also [KPS]). 

\m
We take this opportunity to express our appreciation to  Alberto Elduque for his interest
in our work and to Arun Ram for several helpful discussions.     
 
\newpage   

\head
\S 1. Preliminaries on Quantum Groups \endhead
\b
\subhead {The quantum group $U_q(\g)$} \endsubhead
\m 
Let $\g$ be a finite-dimensional complex simple Lie algebra
with Cartan subalgebra $\h$, root system $\Phi$, and simple roots $\Pi = \{\a_i \mid
i = 1, \dots, r\}\subset \Phi$  relative
to $\h$.  Let $\frak A = (a_{\a,\be})_{\a,\be \in \Pi}$ be the corresponding Cartan
matrix.  Thus, $a_{\a,\be} = 2(\be,\a)/(\a,\a)$, and 
there exist relatively prime positive integers, 

$$d_\a = \frac {(\a,\a)} {2}, \tag 1.1$$

\n so the matrix  $(d_\a a_{\a,\be})$ is symmetric. 
These integers are given explicitly by

$$\vbox{\settabs  
\+ \n \hskip 2  truein & \cr
\+ \n \text{A}$_r$, \text{D}$_r$, \text{E}$_r  \;(r = 6,7,8)$ &
$d_{\a_i} = 1$ for all $1 \leq i \leq r$ \cr 
\+ \n \text{B}$_r$ &
$d_{\a_i} = 2$  for all $1 \leq i \leq r-1$, and $d_{\a_r} = 1$ \cr
\+ \n \text{C}$_r$ &
$d_{\a_i} = 1$  for all  $1 \leq i \leq r-1$,  and $d_{\a_r} = 2$ \cr
\+ \n  F$_4$ & 
$d_{\a_i} = 2, i = 1,2,$  and  $d_{\a_i} = 1, i = 3,4$ \cr
\+ \n G$_2$ & 
$d_{\a_1} = 1$ \text {and} $d_{\a_2} = 3.$ \cr
}\tag 1.2$$   

\b
Let $\K$ be a field of characteristic not 2 or 3. 
Fix $q \in \K$  such that $q$ is not a root of unity, and such that
$q^{\frac {1}{2}} \in \K$.    
For $m,n,d \in \Z_{> 0}$, let

$$[m]_d = \frac{{q^{md} - q^{-md}}} {{q^{d} -
q^{-d}}}, \quad \quad \text {and set} \tag 1.3$$

\n   
$$[m]_d! = \prod_{j = 1}^m [j]_d \tag 1.4$$

\n and  $[0]_d! = 1$.  Let

$$\left [ {{m} \atop {n}} 
\right]_d = \frac{[m]_d!}
{[n]_d! [m-n]_d!}.\tag 1.5$$ 

\n \proclaim{Definition 1.6} The quantum
group $U = U_q(\g)$ is the unital associative algebra over $\K$
generated by elements
$E_\a, F_\a, K_\a$, and $K_\a^{-1}$ (for all $\a \in \Pi$)
 and subject to the relations,
\b
\item {}{(Q1)} $K_\a K_\a^{-1} = 1 = K_\a^{-1}K_\a, \quad  \quad
K_\a K_\be = K_\be K_\a$
\m
\item {}{(Q2)} $K_\a E_\be K_\a^{-1} = q^{(\a,\be)}E_\be$
\m
\item {}{(Q3)} $K_\a F_\be K_\a^{-1}
= q^{-(\a,\be)}F_\be$
\m
\item {}{(Q4)} $E_\a F_\be - F_\be E_\a =
\delta_{\a,\be}\displaystyle{ \frac{ K_\a-K_\a^{-1}}{q^{d_\a}-q^{-d_\a}}}$ 
\m
\item {}{(Q5)}$\displaystyle \sum_{s = 0}^{1-a_{\a,\be}} (-1)^s
\left [{{1-a_{\a,\be}} \atop {s}}
\right]_{d_\a} E_\a^{1-a_{\a,\be}-s} E_\be E_\a^ s = 0$ 
\m
\item {}{(Q6)}$\displaystyle \sum_{s = 0}^{1-a_{\a,\be}}
(-1)^s\left [ {{1-a_{\a,\be}} \atop {s}} 
\right]_{d_\a} F_\a^{1-a_{\a,\be}-s} F_\be F_\a^ s = 0.$  \endproclaim 
\b
Corresponding to each $\l = \sum_{\a \in \Pi} m_\a \a$ in the root lattice $\Z \Phi$ 
there is an element 

$$K_\l = \prod_{\a \in \Pi} K_\a^{m_\a},$$

\n in $U_q(\g)$, and $K_\l K_\mu = K_{\l + \mu}$ for all $\l,\mu \in \Z \Phi$.  
Using relations (Q2),(Q3), we have

$$K_\l E_\be K_\l^{-1} = q^{(\l,\be)}E_\be \quad \quad
\text{and} \quad \quad K_\l F_\be K_\l^{-1} = q^{-(\l,\be)}F_\be \tag 1.7$$

\n for all $\l \in \Z \Phi$ and $\be \in \Pi$.

\b 
The algebra $U_q(\g)$ has a noncommutative Hopf structure 
with comultiplication $\Delta$, antipode $S$, and
counit $\e$  given by

$$\gathered 
\Delta(K_\a) = K_\a \ot K_\a, \quad \quad  \\
\Delta(E_\a) = E_\a \ot 1 + K_\a \ot E_\a 
\\
\Delta(F_\a ) = F_\a  \ot K_\a^{-1}  + 1 \ot F_\a  
\\ 
S(K_\a) = K_\a^{-1}, \\
S(E_\a) = - K_\a^{-1} E_\a  \quad \quad S(F_\a ) =
-  F_\a K_\a \\
\e(K_\a) =1, \quad \quad  \e(E_\a) = \e(F_\a) = 0.   \\
\endgathered \tag 1.8 $$ 
  
\b
The 
opposite comultiplication $\Delta^{\text {op}}$ has the property that if
 
$$\Delta(a) = \sum_a a_{(1)} \ot a_{(2)}, \quad \text{then}
\quad  \Delta^{\text {op}}(a) = \sum_a  a_{(2)} \ot  a_{(1)}.$$

\n (We are adopting
the commonly used Sweedler notation for components of
the comultiplication.)
The algebra $U$ with the same multiplication, same
identity, and  same counit, but with comultiplication given by $\Delta^{\text {op}}$
and with $S^{-1}$ as its antipode is also a Hopf algebra.  The comultiplication  
$\Delta^{\text {op}}$ can be regarded as the composition $\tau \circ \Delta$,
where $\tau: U \ot U \rightarrow U \ot U$ is given by $\tau(a \ot b) = b \ot a$. 
The map $\psi: U \rightarrow U_q(\g)$ with $E_\a \mapsto F_\a$, $F_\a \mapsto E_\a$,
$K_\a \mapsto K_\a^{-1}$ is an algebra automorphism, and a coalgebra anti-automorphism, so
$U$ is isomorphic to
$U_q(\g)$ as a
Hopf algebra (see [CP, p. 211]).  
\b
The defining property of the antipode $S$ in any Hopf algebra is that it is an inverse to
the identity map id$_U$ with respect to the convolution product, so that for all
$a
\in U$,

$$\sum_a S(a_{(1)})a_{(2)} = \e(a) \text{id}_U = \sum_a a_{(1)} S(a_{(2)}). \tag 1.9$$  
\b 
 
\subhead
{Representations}
\endsubhead
\m For any Hopf algebra $U$ the comultiplication operator allows us to define
a
$U$-module structure on the tensor product $V \ot W$ of two $U$-modules  $V$ and $W$, where

$$a (v \ot w) = \sum_a a_{(1)}v \ot a_{(2)}w.$$ 
 
\n The counit $\e$ affords a representation of $U$ on the
one-dimensional module given by $a 1 = \e(a) 1$.
If $V$ is a finite-dimensional $U$-module, then
so is its dual space $V^*$ where the $U$-action is given
by the antipode mapping,

$$(a f)(v) = f( S(a) v),\tag 1.10$$

\n for all $a \in U$, $f \in V^*$, and $v \in V$.  
\b 
A finite-dimensional module $M$ for $U = U_q(\g)$ is the sum of 
weight spaces, $M = \bigoplus_{\nu,\chi} M_{\nu,\chi}$, where
$M_{\nu,\chi} = \{x \in M \mid K_\a x = \chi(\a)q^{(\nu,\a)}x$ for all $\a \in \Pi
\}$, and $\chi: \Z \Phi \rightarrow \{\pm 1\}$ is a group homomorphism.
All the modules considered in this paper will be of 
type I, that is $\chi(\a) = 1$ for all $\a \in \Pi$,
so we will drop $\chi$ from the notation and from our considerations.  
\b 
Each finite-dimensional irreducible $U$-module $M$ has a highest weight
$\l$, which is a dominant integral weight for $\g$ relative
to $\Pi$,  and a unique (up to scalar
multiple) maximal vector
$v^+$ so that
$E_\a v^+ = 0$ and $K_\a v^+ = q^{(\l,\a)}v^+$ for all $\a \in \Pi$.
In particular, the
trivial $U$-module, $\K1,$  has highest weight 0.  We denote by $L(\l)$
the irreducible $U$-module with highest weight $\l$ (and $\chi \equiv 1$.)

\b
When $V$ is a finite-dimensional $U$-module,  
we can suppose $\{x_{\nu,i}\}$ is a basis of $V$ of weight vectors 
 with $i = 1, \dots, \dim V_\nu$.  
Let $\{x_{\nu,i}^*\}$ be the dual basis in $V^*$ so that $x_{\mu,i}^*(x_{\nu,j}) =
\delta_{\mu,\nu}\delta_{i,j}$.   Then 

$$\aligned
(K_\a x_{\mu,i}^*)(x_{\nu,j}) & = x_{\mu,i}^*(S(K_\a)  x_{\nu,j}) =
x_{\mu,i}^*(K_\a^{-1} 
 x_{\nu,j}) \\ & = q^{-(\nu,\a)}x_{\mu,i}^*(x_{\nu,j}) =
q^{-(\nu,\a)}\delta_{\mu,\nu}\delta_{i,j} =  q^{-(\mu,\a)}\delta_{\mu,\nu}\delta_{i,j}, \\
\endaligned \tag 1.11$$

\n from which we see that $x_{\mu,i}^*$ has weight $-\mu$.      
\b 
Suppose $M$ and $N$ are $U$-modules for
a quantum group $U$.  Define a $U$-module
structure on $\Hom_{\K}(M,N)$ by

$$(a  \ga)(m) = \sum_a  a_{(1)} \ga (S(a_{(2)}) m). \tag 1.12$$
 
\n
(This just means that the natural map 
$N \ot M^* \rightarrow \Hom_{\K}(M,N)$ is a $U$-module
homomorphism.)  

Note that if $\ga \in \Hom_U(M,N)$, then

$$(a \ga)(m) = \sum_a a_{(1)} \ga (S(a_{(2)}) m)
= \sum_a a_{(1)} S(a_{(2)})\ga (m) = \e(a) \ga(m)$$

\n so that $a  \ga = \e(a) \ga$.  That says

$$\Hom_U(M,N) \subseteq \Hom_{\K}(M,N)^U = \{ \ga
\in \Hom_{\K}(M,N) \mid a  \ga = \e(a) \ga\}.$$

\n Conversely, if $\ga$ belongs to the invariants
$\Hom_{\K}(M,N)^U$, then we have

$$\aligned
0 = \e(K_\a)\ga& = \ga = K_\a \circ \ga \circ K_\a^{-1}, \quad \text {and hence,}
\\ 0 & = K_\a \circ \ga - \ga \circ K_\a \\
0 = \e(E_\a)\ga & = E_\a  \ga = E_\a \circ \ga 
- K_\a \circ \ga \circ K_\a^{-1} \circ E_\a, \quad  \text {which implies}\\
0 &   = E_\a\circ \ga - \ga \circ E_\a
\\ 0 = \e(F_\a) \ga & =  F_\a  \ga =
F_\a \circ \ga \circ K_\a -  \ga\circ F_\a  \circ  K_\a 
\quad  \text {which implies}\\ 
0 & = F_\a
\circ
\ga - \ga \circ F_\a 
\\
\endaligned $$

\n because of the relations in (1.8).  Since these elements
generate $U_q(\g)$, we see that $\ga \in
\Hom_{\K}(M,N)^U$ implies
$\ga \in \Hom_U(M,N)$, so the two spaces are
equal:

$$\Hom_U(M,N) =  \Hom_{\K}(M,N)^U. \tag 1.13$$ 

Now if $M = N = L(\l)$, a finite-dimensional irreducible $U$-module of highest
weight $\l$, then any  $f \in \Hom_U(M,M)$ must map the highest
weight vector $v^+$ to a multiple of itself.  Since $v^+$ generates $M$
as a $U$-module, it follows that $f$ is a scalar multiple of the identity.
Thus, $\Hom_U(M,M) = {\K}\text{id}_M$.   
So the space of invariants $\Hom_{\K}(M,M)^U$, or equivalently
$(M \ot M^*)^U$, is one-dimensional.  
\b 
\subhead {R-matrices} \endsubhead
\b
The quantum group $U = U_q(\g)$ has a triangular decomposition $U = U^- U^0 U^+$,
where $U^+$ (resp. $U^-$)  is the subalgebra generated by the $E_\a$, 
(resp. $F_\a$), and $U^0$ is the subalgebra generated by the 
$K_\a^{\pm 1}$,  for all $\a \in \Pi$.  Then $U^+ = \sum_\mu U^+_\mu$
where $U^+_\mu = \{a \in U^+ \mid K_\a a K_\a^{-1} = q^{(\mu,\a)}a \}$.
It is easy to see using the automorphism $\psi$ above, that  
$U^-$ has the decomposition $U^- = \sum_{\mu} U^-_{-\mu}$, and
$U^+_\mu \neq (0)$ if and only if $U^-_{-\mu} \neq (0)$.  
There is a nondegenerate bilinear form $(\,,\,)$ on $U$ (see [Jn, Lemma 6.16, Prop. 6.21])
which satisfies 
\m
\n (1.14)
\item{}{(a)} $(b,a) = q^{(2\rho,\mu-\nu)}(a,b)$ for all
$a \in U^-_{-\nu}U^0U^+_{\mu}$ and $b \in U^-_{-\mu}U^0U^+_{\nu}$,
where $\rho$ is the half-sum of the positive roots of $\g$.  
\m
\item{}{(b)} $(F_\a,E_\be) = -\delta_{\a,\be}(q_\a - q_\a^{-1})^{-1}$.   
\b
Choose a basis $a_1^\mu, a_2^\mu, \dots, a_{r(\mu)}^\mu$ of $U_\mu^{+}$
and a corresponding dual basis 
$b_1^\mu, b_2^\mu, \dots, b_{r(\mu)}^\mu$ of $U_{-\mu}^{-}$ so
that $(b_i^\mu,a_j^\mu) = \delta_{i,j}$.  Set

$$ \Theta_\mu = \sum_{i = 1}^{r(\mu)} b_i^\mu \ot a_i^\mu.$$

\n In particular,  $\Theta_{\a} = -(q_\a - q_\a^{-1})F_\a \ot E_\a$ for
all $\a \in \Pi$.  
\b
Suppose $\pi: U_q(\g) \rightarrow \End(V)$ and $\pi' : U_q(\g) \rightarrow
\End(V')$ are two finite-dimensional $U_q(\g)$-modules. Since the set of
weights of $V$ is finite, there are only finitely many $\mu \in \Z \Phi$ which
are differences of weights of $V$.  
Therefore, all but
finitely many $\Theta_\mu$ act as zero on $V$ and $V'$, and 
 $\Theta_{V,V'} = (\pi\ot \pi')(\Theta)$ 
for $\Theta = \sum_{\mu} \Theta_\mu$ is a well-defined mapping
on $V \ot V'$.  By suitably ordering the basis of $V \ot V'$,
we see that each $\Theta_\mu$ for $\mu \neq 0$ acting on $V \ot V'$ has a strictly
upper triangular matrix relative to that basis.  Since $\Theta_0 = 1 \ot 1$, the
transformation $\Theta_{V,V'}$ is unipotent. 
\b
Consider the map $f = f_{V,V'}$, which is defined by  $f_{V,V'}: V \ot V' 
\rightarrow V \ot V'$,  $f(v \ot w) = q^{-(\l,\mu)} v \ot w$ for all $v \in
V_\l$ and  $w \in V'_\mu$.  
\b
\proclaim {Proposition 1.15} (Compare [Jn, Thm. 7.3].)  Let $V,V'$ be finite-dimensional
$U_q(\g)$-modules, and let $\sigma = \sigma_{V',V}$ be defined by

$$\aligned \sigma : V' \ot V & \rightarrow V \ot V'  \\
  w \ot v & \mapsto v \ot w. \\
\endaligned $$

\n  Then the mapping  $\check R_{V',V}: V' \ot V \rightarrow
V \ot V'$  given by
$\check R_{V',V}  = \Theta \circ f \circ \sigma$ is a $U_q(\g)$-module 
isomorphism. 

\endproclaim 

\b
There is an alternate expression for $\Theta_\mu$ which is more explicit and
can be best understood using the
braid group action.  The braid group ${\Cal B}_\g$ associated to $\g$ has generators
$T_\a$, $\a \in \Pi$, and defining relations 

$$(T_\a T_\be)^{m_{\a,\be}} = (T_\be T_\a)^{m_{\a,\be}}$$

\n where $m_{\a,\be} = 2,3,4,6$  if $a_{\a,\be}a_{\be,\a} = 0,1,2,3$ respectively.  
It acts by algebra automorphisms on $U_q(\g)$ according to the following
rules (see [Jn, Sec. 8.14]):

$$\gathered
T_\a(K_\be) = K_{s_\a\be}, \quad \quad T_\a(E_\a) = -F_\a K_\a,
\quad \quad T_\a(F_\a) = -K_\a^{-1}E_\a \\
T_\a (E_\be) = \sum_{t = 0}^{-a_{\a,\be}} (-1)^{t}q_\a^{-t}
(E_\a)^{(-a_{\a,\be}-t)}E_\be (E_\a)^{(t)} \quad {\a \neq \be} \\
T_\a (F_\be) = \sum_{t = 0}^{-a_{\a,\be}} (-1)^{t}q_\a^t 
(F_\a)^{(t)}F_\be (F_\a)^{(-a_{\a,\be}-t)} \quad {\a \neq \be},\\
\endgathered \tag 1.16
$$

\n where $s_{\a}$ is the simple reflection corresponding to $\a$  in the Weyl group $W$ of
$\g$, and

$$(E_\a)^{(n)} = \frac{E_\a^n}{[n]_{d_\a}!},\quad \quad (F_\a)^{(n)} =
\frac{F_\a^n}{[n]_{d_\a}!},$$

\n (the factorials are as defined in (1.4)).  
\b
 Let
$s_i = s_{\a_i}$ be the reflection in the hyperplane perpendicular to the simple root $\a_i$
and fix a reduced decomposition  $w_0 = s_{i_1}s_{i_2} \cdots s_{i_m}$ 
of the longest element of $W$.  Every positive root occurs precisely
once in the sequence

$$\be_1 = \a_{i_1}, \quad \be_2 = s_{i_1}(\a_{i_2}), \; \dots, \;
\be_m = s_{i_1}s_{i_2} \cdots s_{i_{m-1}}(\a_{i_m}). \tag 1.17 $$

The elements

$$
\aligned
E_{\be_t} & = T_{\a_{i_1}}T_{\a_{i_2}} \cdots T_{\a_{i_{t-1}}}(E_{\a_{i_t}}) \\
F_{\be_t} & = T_{\a_{i_1}}T_{\a_{i_2}} \cdots T_{\a_{i_{t-1}}}(F_{\a_{i_t}}) \\
\endaligned \tag 1.18
$$

\n for $t = 1, \dots, m$, 
are called the positive (resp. negative) root vectors of $U_q(\g)$.  
For any $m$-tuple ${\u \ell} = (\ell_1, \dots, \ell_m)$  of
nonnegative integers let

$$\aligned
E^{\u \ell} & = E_{\be_m}^{\ell_m} E_{\be_{m-1}}^{\ell_{m-1}} \cdots
E_{\be_1}^{\ell_1} \\
F^{\u \ell} & = F_{\be_m}^{\ell_m} F_{\be_{m-1}}^{\ell_{m-1}} \cdots
F_{\be_1}^{\ell_1}. \\ 
\endaligned$$

\n Then the elements $E^{\u \ell}$ (resp. $F^{\u \ell}$) determine a basis of $U^+$ (resp.
$U^-$). (See for example, [Jn, Thm. 8.24].)   
\b
Consider for $t = 1, \dots, m$ the sum

$$\Theta^{[t]} = \sum_{j \geq 0} (-1)^j q_\a^{-j(j-1)/2} \frac{(q_\a-q_\a^{-1})^j}{[j]_{d_\a}!}
F_{\be_t}^j \ot E_{\be_t}^j $$

\n where $\a = \a_{i_t}$ and $E_{\be_t}$ and $F_{\be_t}$ are as in (1.18).  
This lies in the direct product of all the $U^-_{-\mu} \ot U^+_{\mu}$.  
Then (as discussed in [Jn, Sec. 8.30]),  $\Theta_\mu$ is the $(U_\mu^- \ot U_\mu^+)$-component
of the product

$$\Theta^{[m]}\Theta^{[m-1]} \cdots \Theta^{[2]}\Theta^{[1]}.$$

\n That component involves only finitely many summands.  
\b
\b
\head  {\S 2.  8-dimensional Representations of $U_q($D$_4)$} \endhead
\b
In this section we
specialize to the case that $\g$ is a finite-dimensional simple complex Lie algebra
of type D$_4$. We can identify $\g$ with the special orthogonal Lie algebra
$so(8)$. The set of roots of $\g$ is given by $\{\pm \e_i \pm \e_j \mid 
1 \leq i \neq j \leq 4\}$, where $\e_i, i = 1, \dots, 4$, is an orthonormal
basis of $\text{\bf R}^4$, and the simple roots may be taken to
be $\a_1 = \e_1 -\e_2, \a_2 = \e_2 - \e_3, \a_3 = \e_3 - \e_4, \a_4 = \e_3 +
\e_4$.    The corresponding
fundamental weights have the following expressions: 

$$\aligned \omega_1 & = \a_1 + \a_2 + (1/2)\a_3 + (1/2)\a_4
=\e_1, \\
\omega_2 & = \a_1 + 2 \a_2 + \a_3 + \a_4 = \e_1 + \e_2 \\
\omega_3 & =
(1/2)\a_1 + \a_2 + \a_3 + (1/2)\a_4 = (1/2)\Big(\e_1+\e_2 + \e_3 - \e_4\Big) \\
\omega_4 & = (1/2)\a_1 + \a_2 + (1/2)\a_3 + \a_4 = (1/2)\Big(\e_1+\e_2 +
\e_3 +
\e_4\Big),
\endaligned \tag 2.1$$
\m
\n and they satisfy $(\omega_i,\a_j) = \delta_{i,j}$ for all $i,j = 1,
\dots, 4$. 
\b
The Lie algebra $\g$ has three 8-dimensional representations - its
natural representation and two spin representations  with highest
weights $\omega_1,\omega_3,$ and $\omega_4$ respectively, and so
does its quantum counterpart $U_q($D$_4)$. ({\it Note we are writing
$\uqd$ rather than $U_q(\g)$ to emphasize the role
that D$_4$ is playing here.})    In particular, for the
natural representation $V$ of $U_q($D$_4)$, there is a basis
$\{x_\mu \mid \mu = \pm \e_i, i = 1, \dots, 4\}$ (compare [Jn, Chap. 5]) such that

$$\gathered K_\a x_\mu = q^{(\mu,\a)} x_\mu, \\
E_\a x_\mu =   
\cases 
0 \ \ \text {\quad if \quad  $(\mu,\a) \neq -1$} \\
x_{\mu + \a} \ \ \text {\quad if \quad  $(\mu,\a) = -1$} \\
\endcases \\
F_\a  x_\mu =  
\cases 
0 \ \ \text {\quad if \quad  $(\mu,\a) \neq 1$} \\
x_{\mu - \a} \ \ \text {\quad if \quad  $(\mu,\a)
= 1$.} \\
\endcases \\
\endgathered \tag 2.2$$

\b
For D$_4$, the group of graph automorphisms of its Dynkin
diagram can be identified with the symmetric
group
$S_3$ of permutations on $\{1,2,3,4\}$ fixing 2. Thus,
$\phi \in S_3$ permutes the simple
roots $\phi(\a_i) = \a_{\phi i}$ and fixes $\a_2$.   Each graph
automorphism
$\phi$ induces an automorphism, which we also denote by $\phi$, on $\uqd$
defined by

$$ E_{\a_i} \mapsto E_{\a_{\phi i}}  \quad \quad
 F_{\a_i} \mapsto F_{\a_{\phi i}} \quad \quad
K_{\a_i} \mapsto K_{\a_{\phi i}}.$$

\n Corresponding to $\phi$ there is new representation, $\pi_\phi :\uqd 
\rightarrow \End(V)$ of $\uqd$ on its natural representation $V$
such that

$$\pi_\phi(a)x = \phi(a)x$$

\n for all $a \in \uqd$ and $x \in V$.  We denote the
$\uqd$-module with this action by $V_\phi$.  Then we have

$$\gathered \pi_\phi(K_\a) x_\mu = K_{\phi \a}x_\mu =  q^{(\mu,\phi\a)}
x_\mu = q^{(\phi^{-1}\mu,\a)} 
x_\mu \\ 
\pi_\phi(E_\a) x_\mu = E_{\phi \a} x_\mu =  
\cases 
0 \ \ \text {\quad if \quad  $(\phi^{-1}\mu,\a) = (\mu,\phi\a) \neq -1$} \\
x_{\mu + \phi \a} \ \ \text {\quad if \quad  $(\phi^{-1}\mu,\a)=
(\mu,\phi\a) = -1$} \\
\endcases \\
\pi_\phi(F_\a) x_\mu = F_{\phi \a} x_\mu =  
\cases 
0 \ \ \text {\quad if \quad  $(\phi^{-1}\mu,\a) = (\mu,\phi\a) \neq 1$} \\
x_{\mu - \phi \a} \ \ \text {\quad if \quad  $(\phi^{-1}\mu,\a)=(\mu,\phi\a)
= 1.$} \\
\endcases \\
\endgathered \tag 2.3$$

\n From these formulas the following is apparent:
\b
\proclaim {Proposition 2.4} $V_\phi$ is an irreducible $\uqd$-module
with highest weight $\phi^{-1}\omega_1$.  
In particular, we have the following highest weights:
\m
\item{}{(i)} $\omega_1$ for $\phi =$ id  or $(3\; 4)$,
\s
\item{}{(ii)} $\omega_3$ for $\phi = (1\;3)$ or $(1\;4\;3)$,
\s
\item{}{(iii)} $\omega_4$ for $\phi = (1\;4)$ or $(1\;3\;4)$.
\endproclaim
\b 
The operators $E_\a, F_\a$ map elements in the basis $\{x_\mu\}$
to other basis elements or 0, and each vector $x_\mu$ is a weight
vector of weight $\phi^{-1}\mu$.  For the natural representation
$V$, the basis elements are $x_{\pm \e_i}, i = 1, \dots, 4$, which
we abbreviate as $x_{\pm i}$.  We can view the module $V_{\phi}$ as displayed
below, where the index $i$ on an arrow indicates that $F_{\a_i}$ maps the higher
vector down to the lower vector and $E_{\a_i}$ maps the lower vector up to the
higher one.  The action of all $F_{\a_j}$'s and $E_{\a_j}$'s not shown is 0.

$$\aligned  & V_\phi \\
& \\
x_1 &\circ \\
& \downarrow \phi^{-1}1\\
x_2 &\circ  \\
& \downarrow 2  \\
x_3&\circ \\
\phi^{-1}3\swarrow &\hskip .2 truein  \searrow \phi^{-1}4\\
x_4\,\circ \quad \quad & \hskip .5 truein  \circ \,x_{-4} \\
\phi^{-1}4\searrow & \hskip .2 truein \swarrow \phi^{-1}3\\
x_{-3}&\circ \\
 & \downarrow 2 \\
x_{-2}&\circ \\
& \downarrow \phi^{-1}1 \\
x_{-1}&\circ \\
\endmatrix  \tag 2.5$$
\m 
\n In particular, taking the identity, $\theta = (1 \;4 \;3)$, and  $\theta^2 =
(1\;3\;4)$  we have the following pictures for the natural representation
$V$, and the spin representations $V_-$ and $V_+$:

$$\matrix  & V & \\
& & \\ & x_1  \circ & \\
 & \downarrow  1 & \\
& x_2  \circ & \\
& \downarrow 2 & \\
& x_3 \circ & \\
3\swarrow &\hskip -.3 truein & \searrow 4\\
x_4\,\circ & \quad \quad & \circ \,x_{-4} \\
4\searrow & &  \swarrow 3\\
& x_{-3} \circ &\\
 & \downarrow 2 &\\
& x_{-2} \circ & \\
& \downarrow  1 &\\
& x_{-1} \circ & \\
\endmatrix \hskip .6 truein
\matrix & V_-& \\ & & \\  & x_1  \circ & \\
 & \downarrow  3 & \\
& x_2  \circ & \\
& \downarrow 2 & \\
& x_3 \circ & \\
4\swarrow &\hskip -.3 truein & \searrow 1\\
x_4\,\circ & \quad \quad & \circ \,x_{-4} \\
1\searrow & &  \swarrow 4\\
& x_{-3} \circ &\\
 & \downarrow 2 &\\
& x_{-2} \circ & \\
& \downarrow  3 &\\
& x_{-1} \circ & \\
\endmatrix \hskip .6 truein  \matrix & V_+& \\ & & \\ & x_1  \circ & \\
 & \downarrow  4 & \\
& x_2  \circ & \\
& \downarrow 2 & \\
& x_3 \circ & \\
1\swarrow &\hskip -.3 truein & \searrow 3\\
x_4\,\circ & \quad \quad & \circ \,x_{-4} \\
3\searrow & &  \swarrow 1\\
& x_{-3} \circ &\\
 & \downarrow 2 &\\
& x_{-2} \circ & \\
& \downarrow  4 &\\
& x_{-1} \circ & \\
\endmatrix  \tag 2.6$$
\m
 
There is a unique (up to scalar multiple) $\uqd$-module homomorphism
$\ast: V_- \ot V_+ \rightarrow V$, which we can compute directly
using the expressions for the comultiplication in (1.8). 
We display the result of this computation 
in Table 2.7 below. 
To make the table more symmetric, we scale the mapping
by $q^{\frac {3}{2}}$.  

$$\table  {\bf\it Table 2.7}
\hfill ! \hfill \hfill ! \hfill \hfill ! \hfill \hfill ! \hfill \hfill
! \hfill \hfill ! \hfill \hfill ! \hfill \hfill ! \hfill  \r 
  \hfill  $\ast$ \hfill ! \hfill $x_1$ \hfill ! 
\hfill $x_2$ \hfill ! \hfill $x_3$ \hfill ! \hfill $x_{-4}$ \hfill ! \hfill $x_4$ \hfill !
\hfill $x_{-3}$ \hfill ! \hfill $x_{-2}$ \hfill ! \hfill
$x_{-1}$  \hfill \r
\hfill ! \hfill \hfill ! \hfill \hfill ! \hfill \hfill ! \hfill \hfill
! \hfill \hfill ! \hfill \hfill ! \hfill \hfill ! \hfill  \rr 
 \hfill  $x_1$ \hfill ! \hfill 0 \hfill ! \hfill 0 
\hfill ! \hfill 0 \hfill ! \hfill $q^{-\frac {3}{2}}x_1$ \hfill ! \hfill 0 \hfill ! \hfill
$q^{-\frac {3}{2}}x_2$ \hfill ! \hfill
$q^{-\frac {3}{2}}x_3$ \hfill ! \hfill $q^{-\frac {3}{2}}x_{-4}$ \hfill \rr
 \hfill   $x_2$ \hfill ! \hfill  0 \hfill ! \hfill 0 \hfill ! \hfill $-q^{-\frac {1}{2}}x_1$
\hfill ! \hfill 0 \hfill ! \hfill \hfil $-q^{-\frac {1}{2}}x_2$  \hfill ! \hfill
$0$ \hfill ! \hfill $q^{-\frac {3}{2}}x_4$ \hfill ! \hfill $q^{-\frac
  {3}{2}}x_{-3}$
\hfill \rr
 \hfill   $x_3$ \hfill ! \hfill 0 \hfill ! \hfill $q^{\frac {1}{2}}x_1$ \hfill ! \hfill 0
\hfill ! \hfill 0 \hfill ! \hfill \hfil $-q^{-\frac {1}{2}}x_3$ \hfill ! \hfill
\hfil $-q^{-\frac {1}{2}}x_4$ \hfill ! \hfill
0 \hfill ! \hfill $q^{-\frac {3}{2}}x_{-2}$  \hfill \rr
 \hfill  $x_4$ \hfill ! \hfill $-q^{\frac {3}{2}}x_1$ \hfill ! \hfill 0 \hfill ! \hfill 0
\hfill ! \hfill 0 \hfill ! \hfill $-q^{-\frac {1}{2}}x_{-4}$\hfil \hfill ! \hfill
$-q^{-\frac {1}{2}}x_{-3}$\hfil \hfill ! \hfill
$-q^{-\frac {1}{2}}x_{-2}$\hfil \hfill ! \hfill   0  \hfill \rr
 \hfill  $x_{-4}$ \hfill ! \hfill 0 \hfill ! \hfill $q^{\frac {1}{2}}x_2$ \hfill ! \hfill
$q^{\frac {1}{2}}x_3$ \hfill ! \hfill $q^{\frac {1}{2}}x_4$ \hfill ! \hfill 0 \hfill ! \hfill
0 \hfill ! \hfill 0 \hfill ! \hfill $q^{-\frac {3}{2}}x_{-1}$  \hfill  \rr
   $x_{-3}$ \hfill ! \hfill $-q^{\frac {3}{2}}x_2$  \hfill ! \hfill 0 \hfill ! \hfill
$q^{\frac {1}{2}}x_{-4}$ \hfill ! \hfill 
$q^{\frac {1}{2}}x_{-3}$ \hfill ! \hfill 0 \hfill ! \hfill
$0$ \hfill ! \hfill  $-q^{-\frac {1}{2}}x_{-1}$ \hfill ! \hfill 0
\hfill  \rr
 \hfill  $x_{-2}$ \hfill ! \hfill  $-q^{\frac {3}{2}}x_3$  \hfill ! \hfill $-q^{\frac
     {3}{2}}x_{-4}$\hfil \hfill ! \hfill 0 \hfill ! \hfill  
$q^{\frac {1}{2}}x_{-2}$ \hfill ! \hfill 0 \hfill ! \hfill
$q^{\frac {1}{2}}x_{-1}$ \hfill ! \hfill
0 \hfill ! \hfill
0  \hfill  \rr
 \hfill   $x_{-1}$ \hfill ! \hfill $-q^{\frac {3}{2}}x_4$ \hfill ! \hfill $-q^{\frac
{3}{2}}x_{-3}$ \hfill ! \hfill $-q^{\frac {3}{2}}x_{-2}$ \hfill ! \hfill 0 \hfill ! \hfill
$-q^{\frac {3}{2}}x_{-1}$ \hfill ! \hfill
0 \hfill ! \hfill
0 \hfill ! \hfill 0 \hfill 
\caption{}
$$

It can be seen readily from (1.8) that the relations

$$\aligned
& (\phi \ot \phi) \Delta(a) = \Delta(\phi(a)) \\
& \phi(S(a)) = S(\phi(a)) \\
\endaligned \tag 2.8 $$

\n hold for all $a \in \uqd$ and all graph automorphisms $\phi$.  
Therefore, by identifying $V_-$ with $V_\theta$ and $V_+$ with
$V_{\theta^2}$ we obtain from 

$$a  (x \ast y) = \sum_a \theta(a_{(1)}) x \ast \theta^2(a_{(2)}) y$$

\n that the following relations hold:

$$
\aligned & \pi_{\theta}(a)  (x \ast y) = \sum_a \theta^2(a_{(1)})  x \ast
a_{(2)}
 y \\ 
& \pi_{\theta^2}(a)  (x \ast y) = \sum_a  a_{(1)}  x \ast 
\theta(a_{(2)}) y. \\
\endaligned \tag 2.9$$

\n They allow us to conclude 
\b
\proclaim {Proposition 2.10} The product $\ast: V_- \otimes V_+ \rightarrow V$
displayed in Table 2.7 above is a $\uqd$-module homomorphism, and it
gives $\uqd$-module homomorphisms,

$$\aligned
& \ast: V_+ \otimes V \rightarrow V_- \\
& \ast:  V \otimes V_- \rightarrow V_+. \\
\endaligned $$

\endproclaim

Suppose now that $\zeta = (1 \;3)$ and $\eta = (1 \; 4)$.  Specializing
(2.5) with $\phi$ equal to these permutations gives:

$$\matrix  & V_\zeta & \\
& & \\ & x_1  \circ & \\
 & \downarrow  3 & \\
& x_2  \circ & \\
& \downarrow 2 & \\
& x_3 \circ & \\
1\swarrow &\hskip -.3 truein & \searrow 4\\
x_4\,\circ & \quad \quad & \circ \,x_{-4} \\
4\searrow & &  \swarrow 1\\
& x_{-3} \circ &\\
 & \downarrow 2 &\\
& x_{-2} \circ & \\
& \downarrow  3 &\\
& x_{-1} \circ & \\
\endmatrix \hskip .7 truein
\matrix & V_\eta& \\ & & \\  & x_1  \circ & \\
 & \downarrow  4 & \\
& x_2  \circ & \\
& \downarrow 2 & \\
& x_3 \circ & \\
3\swarrow &\hskip -.3 truein & \searrow 1\\
x_4\,\circ & \quad \quad & \circ \,x_{-4} \\
1\searrow & &  \swarrow 3\\
& x_{-3} \circ &\\
 & \downarrow 2 &\\
& x_{-2} \circ & \\
& \downarrow  4 &\\
& x_{-1} \circ & \\
\endmatrix \tag 2.11$$

\n As these diagrams indicate, the transformation specified by

$$\aligned \jmath : V & \rightarrow V \\
x_i & \mapsto -x_i \quad \quad i \neq \pm 4 \\
  x_4 & \mapsto -x_{-4} \\
 x_{-4} & \mapsto -x_4 \\
\endaligned
\tag 2.12$$

\n determines explicit $\uqd$-isomorphisms,

$$\jmath: V_- \rightarrow V_\zeta \quad \quad \jmath: V_+ \rightarrow V_\eta.
\tag 2.13$$

\n This together with Proposition 2.10 enables to conclude the following:
\b 
 \proclaim {Proposition 2.14} The map  $\cdot \,: V_\zeta \ot V_\eta \rightarrow V$
given by \; $x \ot y \mapsto x\cdot y = \jmath(x) \ast \jmath(y)$
is a $\uqd$-module homomorphism.  \endproclaim
\b
\subhead Eigenvalues of $\check R_{V,V}$ \endsubhead
\b
Suppose now $V = L(\omega_1)$, the natural representation of $\uqd$. 
The module $V^{\ot 2}$
is completely reducible: $V^{\ot 2} = L(2 \omega_1) \oplus L(\omega_2)
\oplus L(0)$.    By determining the eigenvalues of $\check R_{V,V}$ on
the corresponding maximal vectors, we obtain the eigenvalues of
$\check R_{V,V}$ on the summands.   
\b
First, on the  maximal vector 
$x_1 \ot x_1$ of the summand $L(2 \omega_1)$, 
$\check R_{V,V}(x_1 \ot x_1) = q^{-(\e_1,\e_1)} x_1 \ot x_1 = q^{-1}(x_1 \ot x_1)$.  
The vector $x_2 \ot x_1 -q^{-1}x_1 \ot x_2$ is a maximal vector
for $L(\omega_2)$. Then 

$$\aligned
\check R_{V,V}(x_2 \ot x_1 -q^{-1}x_1 \ot x_2)
& = q^{-(\e_1,\e_2)}\Theta(x_1 \ot x_2 -q^{-1}x_2 \ot x_1) \\
& = x_1 \ot x_2 -q^{-1}x_2 \ot x_1 -(q-q^{-1})x_2 \ot x_1 \\
& = x_1 \ot x_2 - q x_2 \ot x_1 = -q(x_2 \ot x_1 -q^{-1}x_1 \ot x_2), \\
\endaligned $$ 

\n so that $\check R_{V,V}$ acts as $-q$ on $L(\omega_2)$.  We will argue 
in the remark following Proposition 2.19 below that $\check R_{V,V}$ acts on $L(0)$
as $q^7$.   Thus, the eigenvalues of $\check R_{V,V}$ are given by

$$\aligned q^{-1}  &  \hskip .5 truein L(2\omega_1) \\
-q &  \hskip .5 truein L(\omega_2) \\
q^7 & \hskip .5 truein L(0), \\ 
\endaligned \tag 2.15
$$
\n (compare [D], [R]). Setting ${\Cal S}^2(V) = L(2\omega_1) \oplus L(0)$ we see that
$\check R_{V,V} + q$id$_V$ maps $V^{\ot 2}$ onto ${\Cal S}^2(V)$, the {\it quantum
symmetric space}.   There are analogous decompositions for the
two spin representations of $\uqd$ obtained by replacing $\omega_1$
by $\omega_3$ and $\omega_4$.  
\b 
For each graph
automorphism
$\phi$, we consider $\check R = \check R_{V_\phi,V_\phi}$ 
 on $V_\phi \ot V_\phi$. Note the bilinear form $(\,,\,)$ on
$\uqd$ in [Jn, (6.12)] is
$(\phi \ot \phi)$-invariant, so

$$\Theta_{\phi(\mu)} = \sum_i \phi(b_i^\mu) \ot \phi(a_i^\mu) = (\phi \ot \phi)(\Theta_\mu).$$

\n Now

$$
\aligned
\check R(v \ot w)&  = (\Theta \circ f \circ \sigma)(v \ot w)
= q^{-(\phi^{-1}\l,\phi^{-1}\nu)} \Theta (w \ot v) \\
& = q^{-(\l,\nu)} \Theta (w \ot v) = \check R_{V,V}(v \ot w), \\
\endaligned \tag 2.16$$

\n whenever $w$ is in $V_\l$ and $v
 \in V_\nu$. 
Consequently, the actions of the R-matrices agree. It follows that

$${\Cal S}^2(V) = {\Cal S}^2(V_\phi) \tag 2.17$$

\n for all graph automorphisms $\phi$.  
\b
\subhead A ``test'' element of ${\Cal S}^2(V)$  \endsubhead 
\b
In the next two chapters we will compute various $\uqd$-module homomorphisms and
prove that they are equal by evaluating them on a particular element
of ${\Cal S}^2(V)$.  A good
``test'' element for these calculations is $x_1 \ot x_{-1} + x_{-1} \ot x_1$,
but first we need to demonstrate that it does in fact belong to ${\Cal S}^2(V)$.  
Since $\check R_{V,V} + q$id$_V$ maps $V^{\ot 2}$ onto ${\Cal S}^2(V)$
for $V = L(\omega_1)$, the natural representation,  
it suffices to show that the image of $x_1 \ot x_{-1}$ under 
$\check R_{V,V} + q$id$_V$ is $x_1 \ot x_{-1} + x_{-1} \ot x_1$. 
Using the fact that $E_\gamma x_1 = 0$ for all positive roots $\gamma$
we obtain

$$\aligned \check R_{V,V}(x_1 \ot x_{-1}) & = (\Theta \circ f \circ \sigma) 
(x_1 \ot x_{-1}) \\ & = q \Theta (x_{-1} \ot x_1)  
= q (x_{-1} \ot x_1).  \\
\endaligned  $$
 
\n Consequently,  
$$(\check R_{V,V} + q\text{id}_V)(x_1 \ot x_{-1}) = q x_{-1} \ot x_1 +
q x_1 \ot x_{-1},$$ 

\n which implies the desired conclusion 

$$x_{-1} \ot x_1 + x_1 \ot x_{-1} \in {\Cal S}^2(V).$$
   
\b
\subhead Bilinear Forms \endsubhead
\b 
For the natural representation $V$ of $U = \uqd$, the dual space $V^*$ has
weights $\pm \e_i$, $i = 1, \dots, 4$ by (1.11).   The 
vector $x_{-1}^*$ is a maximal vector since
it corresponds to the unique dominant weight,  $\e_1$,
and   $V^*$ is
isomorphic to $V$.   We have seen in (1.13) that the
space of invariants $(V \ot V^*)^U$ is one-dimensional.  Since $V \cong V^*$,
we have that the space of invariants $(V \ot V)^U$ is one-dimensional.   
Thus, in $V \ot V$ there is a unique copy of the trivial module, and
hence a unique (up to scalar multiple) $\uqd$-module homomorphism

$$(\,|\,): V \ot V \rightarrow \K.$$

\n The condition of $(\,|\,)$ being a $\uqd$-module homomorphism is
equivalent to  

$$(a  x | y) = (x | S(a)  y)  \tag 2.18 $$

\n for all $a \in U$ and $x,y \in V$ (see [Jn, p.122]).  There is an analogous
bilinear form giving the unique (up to multiples) $\uqd$-module
homomorphism, $V_\phi \ot V_\phi
\rightarrow \K$. All these bilinear forms can be
regarded as $\uqd$-module maps from the symmetric
tensors ${\Cal S}^2(V_\phi) \rightarrow \K$, since the trivial module
is an irreducible summand of ${\Cal S}^2(V_\phi)$ in each case.  

Now for $a \in \uqd$, 

$$
\aligned \sum_a (\pi_\phi(a_{(1)}) x | \pi_\phi(a_{(2)}) y) & =
\sum_a (\phi(a_{(1)}) x | \phi(a_{(2)}) y) \\
& = \sum_a (x | S(\phi(a_{(1)}))\phi(a_{(2)}) y) \\
& = \sum_a (x | \phi( S(a_{(1)})a_{(2)})  y \\
& = (x | \phi \Big(\sum_a S(a_{(1)})a_{(2)} \Big) ) y \\
& = (x | \phi (\e(a)1)  y) = \e(a)(x | y). \\
\endaligned \tag 2.19$$
 
\n (We have used in the third line the
fact that $\phi$ and $S$ commute (2.8).)  Hence, by the calculation in (2.19), and
by direct computation of $(\,|\,)$ on $V \otimes V$,  we have
\b
\proclaim {Proposition 2.20} The $\uqd$-module
homomorphism: $(\,|\,): V \ot V \rightarrow
\K$  determines a $\uqd$-module
homomorphism: $(\,|\,): V_\phi \ot V_\phi \rightarrow
\K$ for any graph automorphism $\phi$.  It is
given by

$$(x_i|x_j) = 
\cases 
0 \ \ \text {\quad if \quad  $j \neq -i$} \\
(-1)^{i+1}\displaystyle{\frac {q^{-(4-i)}}{q^3 + q^{-3}}} \ \ \text {\quad if \quad  $i
> 0$ and $
j = -i$} \\
\\
(-1)^{i+1}\displaystyle{ \frac {q^{i+4}}{q^3 + q^{-3}}} \ \ \text {\quad if \quad 
$i < 0$ and $
j = -i$}. \\
\endcases$$

\endproclaim
\b
We are free to scale the bilinear form $(\,|\,)$ on
$V_\phi$  by any nonzero element in
$\K$. In writing down the values of $(x_i|x_{-i})$ in
Proposition 2.20, we have chosen a particular
scaling with an eye towards results (such as Proposition 4.2) to follow. 
It is easy to see using the basis elements $x_i$ that the
following holds.

\b
\proclaim {Proposition 2.21} The transformation $\jmath$ is an isometry, 
(\,$(x|y) = (\jmath(x)|\jmath(y))$ for all $x,y \in V$), relative to the bilinear form
in Proposition 2.20. \endproclaim
\b
{\bf Remark.} Observe that $\check R_{V,V} \circ (\,|\,) = \xi (\,|\,)$, where
$\xi$ is the eigenvalue of $\check R_{V,V}$ on $L(0)$.   Since
$\check R_{V,V}(x_1 \ot x_{-1}) = q(x_{-1} \ot x_1)$, as our calculations
for the ``test element'' show, we have $\xi(x_1|x_{-1}) = q(x_{-1}|x_1)$. Therefore, 
$\xi q^{-3} = q q^3$ or $\xi = q^7$, as asserted in (2.15).     
\b
For D$_4$, the half-sum of the positive roots is given by  $\rho =
3\a_1 +5 \a_2 +3\a_3 + 3\a_4$, so $K_{2\rho}^{-1} =
K_{\a_1}^{-6}K_{\a_2}^{-10}K_{\a_3}^{-6}K_{\a_4}^{-6}$. Therefore,
on the natural representation
$V$, $K_{2\rho}^{-1}$ has a diagonal matrix relative
to the basis $\{x_1,x_2,x_3,x_4,x_{-4},x_{-3},x_{-2},x_{-1}\}$
with corresponding diagonal entries given by 
$\{q^{-6},q^{-4},q^{-2},1,1,q^2,q^4,q^6\}$.  It is easy to verify
directly that 

$$(x | y) = (y | K_{2\rho}^{-1} x). \tag 2.22$$ 

\vskip .5 truein

 \head {\S 3. Octonions and Quantum Octonions} \endhead
\b
Suppose
$\zeta = (1\;3)$, $\eta = (1\;4)$, and $V_\phi$ is the representation of $\uqd$
coming from the graph automorphism $\phi$. Then

$$\dim_\K\Hom_{\uqd}({\Cal S}^2(V_\zeta) \ot {\Cal S}^2(V_\eta), \K) = 1, \tag 3.1$$

\n and we can realize these homomorphisms explicitly as multiples of the mapping

$$\aligned
{\Cal S}^2(V_\zeta)&\ot {\Cal S}^2(V_\eta)\subset V_\zeta \ot V_\zeta \ot V_\eta \ot V_\eta
\rightarrow
\K
\\ & \sum_{k,\ell} u_k \ot v_k \ot w_\ell \ot y_\ell \mapsto
\sum_{k,\ell}(u_k|v_k)(w_\ell|y_\ell).
\\
\endaligned \tag 3.2$$

\n On the other hand, the mapping

$$\align {\Cal S}^2(V_\zeta) \ot {\Cal S}^2(V_\eta) &\subset V_\zeta \ot V_\zeta \ot V_\eta \ot
V_\eta
\rightarrow V_\zeta \ot V_\eta \ot V_\zeta \ot V_\eta \rightarrow V \ot V \rightarrow
\K  \tag 3.3 \\ 
\sum_{k,\ell} u_k \ot v_k \ot w_\ell \ot y_\ell  & \mapsto \sum_{k,\ell} u_k \ot \check
R_{V_{\zeta},V_{\eta}}(v_k \ot w_\ell)
\ot y_\ell\\  
& \quad \quad  \mapsto
\sum_{k,\ell} \sum_i u_k\cdot w_\ell^{(i)}
\ot v_k^{(i)} \cdot y_\ell
\mapsto \sum_{k,\ell}\sum_i (u_k\cdot w_\ell^{(i)} | v_k^{(i)} \cdot y_\ell),\\
\endalign $$

\n is a
$\uqd$-module homomorphism, (where ``$\cdot$'' is the product in Proposition 2.14 and
$\check R_{V_{\zeta},V_{\eta}}(v_k \ot w_\ell) = \sum_i w_\ell^{(i)} \ot v_k^{(i)}$
is the braiding morphism).     Thus, it must be a multiple of the mapping 
in (3.2).  

We now consider what would happen if $q = 1$ and compute the multiple in that
case.  Here the braiding morphism may be taken to be
$v \ot w \mapsto w \ot v$ since that is a module map in the $q = 1$ case. 
We evaluate the
mappings at $e \ot e \ot e \ot e$ where $e = x_4-x_{-4}$.  Now by Table 2.7, we have

$$\aligned
e^2 & = \jmath (e) \ast \jmath (e) = (x_4 - x_{-4}) \ast 
(x_4 - x_{-4}) \\ & = q^{\frac{1}{2}}x_4 - q^{-\frac{1}{2}}x_{-4}, \\
\endaligned \tag 3.4$$ 

\n which shows that at $q = 1$ that $e^2 = e$.   By Proposition 2.19, 

$$(e|e) = -(x_4|x_{-4}) - (x_{-4}|x_4) = \frac {2}{q^3 + q^{-3}}, \tag 3.5$$

\n which gives $(e|e)_{{}_1} = 1$ at $q = 1$.
{\it (Here and throughout the paper we use
$(\,|\,)_{{}_1}$ to denote
the bilinear form at $q = 1$.)}    Therefore, $(e|e)_{{}_1}(e|e)
_{{}_1} = 1$ and
$(e^2|e^2)_{{}_1} = (e|e)_{{}_1} = 1$.  Thus, the two maps in (3.2)
and (3.3) are equal at $q = 1$.  By
specializing $\sum_k u_k \ot v_k = x \ot x$ and  $\sum_\ell w_\ell \ot y_\ell = y \ot y$
(which are symmetric elements at $q = 1$),  we
have as a consequence,

$$(x\cdot y|x \cdot y)_{{}_1} = (x|x)_{{}_1}(y|y)_{{}_1}.\tag 3.6$$

\n Using the multiplication table (Table 2.7) for $x \ast y$ and the
formula $x\cdot y = \jmath (x) \ast \jmath (y)$, we can verify that
$e = x_4 - x_{-4}$ is the unit element relative to the product ``$\cdot$''
when $q = 1$. (For further discussion, see the beginning of Section 5.)
Thus, the product
$\cdot: V_\zeta \ot V_\eta \rightarrow V$  for $q = 1$ gives an
8-dimensional unital
composition algebra over $\K$ whose underlying vector space is  $V = V_\zeta =
V_\eta$.   Such an algebra must be the octonions. Thus, we have
\b
\proclaim {Proposition 3.7} Let $V$ be the 8-dimensional vector space over
$\K$ with basis $\{x_{\pm i} \mid i = 1, \dots, 4\}$.  Then $V$
with the product $\cdot$ defined  by $x \cdot y = \jmath(x) \ast \jmath(y),$
where $\jmath$ is as in (2.12) and $\ast$ is as in Table 2.7, is the
algebra of octonions when $q = 1$.  \endproclaim 
\b
Because of this result we are motivated to make the following definition.
\b
\subhead {Quantum Octonions} \endsubhead 
\m
\proclaim {Definition 3.8} Let $V$ be the 8-dimensional vector space over
$\K$ with basis $\{x_{\pm i} \mid i = 1, \dots, 4\}$.  Then $V$
with the product ``$\cdot$'' defined  by 

$$x \cdot y = \jmath(x) \ast \jmath(y),$$

\n where $\jmath$ is as in (2.12) and $\ast$ is given by Table 2.7, is
the algebra of {\it quantum octonions}.  We denote this
algebra by  $\Bbb O_q = (V,\cdot)$.  
 \endproclaim
\b 

The multiplication table for the quantum octonions is obtained from
Table 2.7 by interchanging the entries in the columns and rows labelled by $x_{-4}$ and 
$x_4$.   As a result, we have

\b   
$$\table  {\bf \it Table 3.9} 
  \hfill ! \hfill \hfill ! \hfill \hfill ! \hfill
 \hfill ! \hfill \hfill ! \hfill \hfill ! \hfill \hfill ! \hfill \hfill ! \hfill\r
 \hfill  $\cdot$ \hfill ! \hfill $x_1$ \hfill ! \hfill   $x_2$  \hfill ! \hfill  $x_3$   \hfill ! \hfill  $x_{-4}$   \hfill ! \hfill  $x_4$ \hfill ! \hfill
   $x_{-3}$  \hfill ! \hfill  $x_{-2}$ \hfill ! \hfill
$x_{-1}$ \hfill
\r
\hfill ! \hfill \hfill ! \hfill \hfill ! \hfill \hfill !
 \hfill \hfill ! \hfill \hfill ! \hfill \hfill ! \hfill \hfill ! \hfill\rr
  \hfill   $x_1$ \hfill ! \hfill 0 \hfill ! \hfill 0 
\hfill ! \hfill 0 \hfill ! \hfill 0 \hfill ! 
\hfill  $q^{-\frac {3}{2}}x_1$  \hfill ! \hfill $q^{-\frac {3}{2}}x_2$ \hfill ! \hfill
$q^{-\frac {3}{2}}x_3$ \hfill ! \hfill $q^{-\frac {3}{2}}x_{-4}$  \hfill
\rr
  \hfill  $x_2$ \hfill ! \hfill 0 \hfill ! \hfill 0 \hfill ! \hfill $-q^{-\frac {1}{2}}x_1$
\hfill ! \hfill $-q^{-\frac {1}{2}}x_2$   \hfill ! \hfill 0 \hfill ! \hfill
$0$ \hfill ! \hfill $q^{-\frac {3}{2}}x_4$ \hfill ! \hfill $q^{-\frac
  {3}{2}}x_{-3}$ \hfill \rr
  \hfill  $x_3$ \hfill ! \hfill 0 \hfill ! \hfill $q^{\frac {1}{2}}x_1$ \hfill ! \hfill 0
\hfill ! \hfill $-q^{-\frac {1}{2}}x_3$ \hfill ! \hfill 0 \hfill ! \hfill
 $-q^{-\frac {1}{2}}x_4$  \hfill ! \hfill
0 \hfill ! \hfill $q^{-\frac {3}{2}}x_{-2}$  \hfill  \rr
 \hfill   $x_4$ \hfill ! \hfill  0 \hfill ! \hfill  $q^{\frac {1}{2}}x_2$ 
\hfill ! \hfill $q^{\frac {1}{2}}x_3$ \hfill ! \hfill 0 \hfill ! \hfill $q^{\frac
{1}{2}}x_{4}$  \hfill ! \hfill
0 \hfill ! \hfill 0 \hfill ! \hfill $q^{-\frac
{3}{2}}x_{-1}$  \hfill  \rr
 \hfill  $x_{-4}$ \hfill ! \hfill $-q^{\frac {3}{2}}x_1$ \hfill ! \hfill 0 \hfill ! \hfill 0
\hfill ! \hfill
$-q^{-\frac {1}{2}}x_{-4}$ \hfill ! \hfill 0 \hfill ! \hfill
$-q^{-\frac {1}{2}}x_{-3}$ \hfill ! \hfill $-q^{-\frac {1}{2}}x_{-2}$
\hfill ! \hfill 0  \hfill  \rr
  \hfill  $x_{-3}$ \hfill ! \hfill $-q^{\frac {3}{2}}x_2$  \hfill ! \hfill 0 \hfill ! \hfill
$q^{\frac {1}{2}}x_{-4}$ \hfill ! \hfill  0 \hfill ! \hfill $q^{\frac {1}{2}}x_{-3}$  \hfill
! \hfill
$0$ \hfill ! \hfill $-q^{-\frac {1}{2}}x_{-1}$ \hfill ! \hfill 0
\hfill \rr
  \hfill  $x_{-2}$ \hfill ! \hfill $-q^{\frac {3}{2}}x_3$ \hfill ! \hfill  $-q^{\frac
{3}{2}}x_{-4}$ \hfill ! \hfill
    0 \hfill ! \hfill  
0 \hfill ! \hfill $q^{\frac {1}{2}}x_{-2}$ \hfill ! \hfill 
$q^{\frac {1}{2}}x_{-1}$ \hfill ! \hfill 0 \hfill ! \hfill 0  \hfill
\rr
  \hfill  $x_{-1}$ \hfill ! \hfill $-q^{\frac {3}{2}}x_4$ \hfill ! \hfill $-q^{\frac
{3}{2}}x_{-3}$ \hfill ! \hfill $-q^{\frac {3}{2}}x_{-2}$ \hfill ! \hfill 
$-q^{\frac {3}{2}}x_{-1}$ \hfill ! \hfill 0 \hfill ! \hfill 0  \hfill
! \hfill 0 \hfill ! \hfill 0 \hfill 
\caption{} $$

\b 
The product map: $p : V_\zeta \ot V_\eta \rightarrow
V$, $p(x \ot y) = x \cdot y$, of the quantum octonion algebra is the composition
$p =
\ast \circ (j
\ot j)$  of two $\uqd$-module maps, so is itself a $\uqd$-module homomorphism. 
\b
Let  $\pi: \uqd \rightarrow
\End(V)$ be the natural representation of  $\uqd$.  Then we have

$$ \aligned  \pi(a)(x \cdot y) & = a \cdot p(x \ot y) 
= p  ( \Delta(a)(x \ot y)) \\ 
& = p \Big ( \sum_a \zeta (a_{(1)})
x \ot  \eta (a_{(2)})y \Big) \\
& = \sum_a \jmath(\zeta(a_{(1)}x)) \ast \jmath(\eta(a_{(2)}y)) \\
& = \sum_a \Big(\zeta(a_{(1)}) x\Big) \cdot \Big( \eta (a_{(2)}) y\Big). \\
\endaligned \tag 3.10$$

\n   Now when $q = 1$ we can identify
$\uqd$ with the universal enveloping algebra $U(\g)$ of a simple Lie algebra
$\g$ of type
D$_4$.  The comultiplication specializes in this case to the usual comultiplication
on $U(\g)$, which has $\Delta(a) = a \ot 1 + 1 \ot a$ for all $a \in \g$. In
particular, for elements $a \in \g$, (3.10) becomes

$$\pi(a)(x \cdot y) = \big(\zeta(a)x\big) \cdot y + x \cdot 
\big(\eta(a)y\big). 
\tag
3.11$$

\n   Equation (3.11) is commonly referred to
as the {\it Principle of Local Triality}  (see for example, [J, p. 8] or
[S, p. 88]).  As a result we
have,    
\b
\proclaim{\bf Proposition 3.12} The quantum octonion algebra $\Bbb O_q = (V, \cdot)$
satisfies the $q$-Principle of Local Triality,

$$\pi(a)(x \cdot y) = \sum_a \Big(\zeta(a_{(1)}) x\Big) \cdot \Big(\eta(a_{(2)}) y\Big)$$

\n for all $a \in \uqd$ and $x,y \in V$, where $\zeta$ and $\eta$ 
are the graph automorphisms $\zeta = (1 \; 3)$ and 
$\eta = (1 \; 4)$ of $\uqd$, and $\Delta(a) = \sum_a a_{(1)} \ot a_{(2)}$
is the comultiplication in $\uqd$.      
\endproclaim 
\b
\subhead {A quantum para-Hurwitz algebra} \endsubhead 

\b At $q =1$, the map $\jmath: V \rightarrow V$ sending $x_i \mapsto -x_{-i}$
for $i \neq \pm 4$,  $x_4 \mapsto -x_{-4}$, and $x_{-4} \mapsto -x_4$,
agrees with the standard involution $x \rightarrow \o x$ on the
octonions so that

$$x \ast y = \jmath(x) \cdot \jmath(y) = \o x \cdot \o y. \tag 3.13$$

\n Thus, the product ``$\ast$'' displayed in Table 2.7 gives what is called the
{\it para-Hurwitz algebra} when $q = 1$.  (Results on
para-Hurwitz algebras can be found in [EM], [EP], [OO1], [OO2].)    
\b
\proclaim {Definition 3.14} We say that the algebra $\Bbb P_q = (V,\ast)$ with the
multiplication $\ast$ given by Table 2.7 is the {\it
quantum para-Hurwitz algebra}.  It satisfies 

$$x \ast y =  \jmath(x) \cdot \jmath(y)$$

\n where $\Bbb O_q = (V,\cdot)$ is the quantum octonion algebra and $\jmath$
is the transformation on $V$ defined in (2.12).   \endproclaim
\b
Our definitions have been predicated on using the mapping
$\cdot: V_\zeta \ot V_\eta \rightarrow V$.   We can ask what would
have happened if we had used the projection $V_\eta \ot V_\zeta \rightarrow V$
instead?  To answer this question, it is helpful to observe that
the mapping $\jmath$ determines a $\uqd$-module isomorphism
$\jmath: V_{\zeta \eta \zeta} = V_{(3 \; 4)} \rightarrow V$, which is
readily apparent from comparing the following picture of $V_{(3 \; 4)}$
with the diagram for $V$ in (2.6).

$$\matrix  & V_{(3 \; 4)} & \\
& & \\ & x_1  \circ & \\
 & \downarrow  1 & \\
& x_2  \circ & \\
& \downarrow 2 & \\
& x_3 \circ & \\
4\swarrow &\hskip -.3 truein & \searrow 3\\
x_4\,\circ & \quad \quad & \circ \,x_{-4} \\
3\searrow & &  \swarrow 4\\
& x_{-3} \circ &\\
 & \downarrow 2 &\\
& x_{-2} \circ & \\
& \downarrow  1 &\\
& x_{-1} \circ & \\
\endmatrix$$ 

Suppose as before that $p: V_\zeta \ot V_\eta \rightarrow V$ is the product
$p(x \ot y) = x \cdot y$ in the quantum octonion algebra.  Consider the
mapping $\jmath \circ p \circ (\jmath \ot \jmath): V_\eta \ot V_\zeta \rightarrow
V$. We demonstrate that this is a $\uqd$-module map.  Bear in mind in these
calculations that $\jmath: V_\zeta \rightarrow V_- = V_{\zeta \eta}$
and $\jmath: V_\eta
\rightarrow V_+ = V_{\eta \zeta}$ for $\zeta = (1\;3)$ and
$\eta = (1\;4)$.  Then for $a \in
\uqd$ we have

$$\aligned 
\jmath \circ p \circ (\jmath \ot \jmath) (a(x \ot y)) & =
\jmath \circ p \circ (\jmath \ot \jmath) \Big( \sum_a \eta(a_{(1)})x \ot \zeta(a_{(2)})y
\Big)\\ & =
\jmath\circ  p \Big( \sum_a  j(\eta(a_{(1)})x) \ot j(\zeta(a_{(2)})y) \Big) \\
& =
\jmath \circ p \Big(\sum_a \eta\zeta(a_{(1)})j(x) \ot  \zeta\eta(a_{(2)})j(y) \Big)\\
& =
\jmath \circ p \Big(\sum_a \zeta^2\eta\zeta(a_{(1)})j(x) \ot 
\eta^2\zeta\eta(a_{(2)})j(y)
\Big)\\ 
\endaligned
$$

\n Now using the fact that $\eta\zeta \eta = \zeta \eta \zeta$ and
that $(\phi \ot \phi) \Delta = \Delta \phi$ for all graph automorphisms
$\phi$ together with (3.10), we see that the last expression can be rewritten as 

$$\aligned
& \jmath \circ p \Big((\zeta \ot \eta)\big(\Delta(\zeta \eta \zeta(a))\big)
\big (\jmath(x)
\ot
\jmath(y)\big) \Big) \\
& \hskip 1 truein = \jmath \Big (\zeta\eta \zeta(a)p\big(\jmath(x)\ot
\jmath(y)\big )\Big)
\quad (p \; \text {is\;a}\;\uqd-\text{homomorphism})\\
& \hskip 1 truein =  a  (\jmath \circ p \circ (\jmath \ot \jmath))(x \ot y). \\
\endaligned $$
\b
Thus, $\jmath \circ p \circ (\jmath \ot \jmath) : V_\eta \ot V_\zeta \rightarrow V$
is a $\uqd$-module homomorphism.   Since the space 
$\Hom_{\uqd}(V_\eta \ot V_\zeta,V)$ is one-dimensional, we may assume the
product $p': V_\eta \ot V_\zeta \rightarrow V$ is given by the map we have just
found, so that $p' = \jmath \circ p \circ (\jmath \ot \jmath)$.  
Alternately, $\jmath \circ p' = p \circ (\jmath \ot \jmath)$.  This
says that had we defined the quantum octonions using the product $p'$ rather than
$p$ so that $x \cdot' y = p'(x \ot y)$, the two algebras would be isomorphic
via the map
$\jmath$. 
\b
\head {\S 4. Properties of Quantum Octonions} \endhead
\b
We use $\uqd$-module homomorphisms to deduce various properties of the quantum
octonions. At $q = 1$ these properties specialize to
well-known identities satisfied by the octonions, which can be found 
for example in [ZSSS].  Since the vector spaces $V$, $V_\zeta$, and
$V_\eta$ are the same, we will denote the identity map on them simply by
``id'' in what follows. 
\b
\proclaim { Proposition 4.1}  The following diagram commutes

$$\matrix  &  (\,|\,) \ot \text{id} & \\
{\Cal S}^2(V_\zeta) \ot V_\eta  & \longrightarrow & V_\eta\\
&  & \\
\text{id} \ot p \downarrow &  & \downarrow \text {id} \\
 & p \circ (\jmath \ot \text{id}) & \\
V_\zeta \ot V & \longrightarrow & V_\eta, \\
\endmatrix $$

\n so that 

$$\sum_k \jmath(u_k)\cdot \big (v_k \cdot y \big) = \sum_k(u_k|v_k)y \tag 4.2$$ 

\n for all $\sum_k u_k \ot v_k \in {\Cal S}^2(\Bbb O_q),\; y \in \Bbb O_q$.  
In particular, at
$q = 1$ this reduces to the well-known identity satisfied by the octonion algebra $\Bbb O$:

$$\o x \cdot \big( x \cdot y\big) = (x | x)_{{}_1}y. \tag 4.3$$

\endproclaim 
\b
{\bf Proof.}  Since $p: V_\zeta \ot V_\eta \rightarrow V$ is a $\uqd$-module 
homomorphism, we have that 
$\pi(a) \circ p = p\Big((\zeta \ot \eta)\Delta(a)\Big)$, which implies that
$\pi(\eta(a)) \circ p = p\Big((\zeta \ot \eta)\Delta(\eta a)\Big)
= p\Big((\zeta \eta \ot \text{id})\Delta(a)\Big)$.  This says that 
$p: V_{\zeta \eta} \ot V \rightarrow V_\eta$ is also a $\uqd$-module
homomorphism.  Since  $\jmath: V_\zeta \rightarrow V_{\zeta \eta}$
is a $\uqd$-module homomorphism, it follows that $\big( p \circ (\jmath \ot \text{id}))
\circ (\text{id} \ot p) \in \Hom_{\uqd}({\Cal S}^2(V_\zeta) \ot V_\eta, \K)$.  As that
space of homomorphisms has dimension one, the maps in the diagram must
be proportional. 

To find the proportionality constant, let us consider the actions using
the test element:
 
$$(x_1 \ot x_{-1} + x_{-1} \ot x_1) \ot x_1 \fimap x_1$$

\n and

$$\aligned
p \circ (\jmath \ot \text{id}) \circ (\text {id} \ot p) & \Big((x_1 \ot x_{-1} + x_{-1} \ot
x_1)\ot x_1
\Big)  =  \\
& p \circ (\jmath \ot \text{id})\Big( -q^{\frac {3}{2}} x_1 \ot x_4 \Big) = q^{\frac
{3}{2}} x_1 \cdot x_4 = x_1.\\
\endaligned $$

\n Therefore, the two maps are equal as claimed.  At $q = 1$
we may substitute  $x \ot x \in {\Cal S}^2(\Bbb O)$ and $y \in \Bbb O$ into the above
identity to obtain the corresponding one for the octonions.  \qed 
\b
\proclaim { Proposition 4.4}  The following diagram commutes

$$\matrix  & \text{id}\ot  (\,|\,) & \\
V_\zeta \ot {\Cal S}^2(V_\eta)  & \longrightarrow & V_\zeta\\
&  & \\
p \ot \text{id}\downarrow &  & \downarrow \text {id} \\
 & p \circ (\text{id}\ot \jmath) & \\
V \ot V_\eta & \longrightarrow & V_\zeta, \\
\endmatrix $$

\n so that 

$$\sum_k \big (y\cdot u_k \big) \cdot \jmath(v_k) = \sum_k(u_k|v_k)y \tag 4.5$$ 

\n for all $\sum_k u_k \ot v_k \in {\Cal S}^2(\Bbb O_q), \; y\in \Bbb O_q$.  In
particular, at
$q = 1$ this gives the identity 

$$\big( y \cdot x\big)\cdot \o x = (x | x)_{{}_1}y \tag 4.6$$

\n which is satisfied by the octonions.

\endproclaim 
\b
{\bf Proof.} Both maps are $\uqd$-module
homomorphisms and so must be multiples of each other. We compute
their images on $x = x_1 \ot (x_1 \ot x_{-1} + x_{-1} \ot x_1)$:

$$\aligned
&\big( \text{id} \ot (\,|\,)\big) x = \Big((x_1 | x_{-1}) + (x_{-1}|x_1)\Big)x_1
= x_1, \\
& p \circ (\text{id} \ot \jmath) \circ (p \ot \text{id})x 
= p\circ (\text{id} \ot \jmath)(q^{-\frac{3}{2}}x_{-4} \ot x_1)
= -q^{-\frac{3}{2}}x_{-4} \cdot x_1 = x_1. \\
\endaligned $$

\n Hence we have the desired equality.  \qed 
\b
\proclaim{Proposition 4.7}  The following diagram commutes

$$\matrix  & & p \ot \text{id} &   \\
& V_\zeta \ot V_\eta \ot V  & \quad \longrightarrow \quad & V \ot V  \\
&  & &  \\
\text{id}\ot\Big(\jmath \circ p \circ (\text{id} \ot \jmath)\Big) & 
\downarrow & \quad
&
\downarrow (\,|\,)\\ & & (\,|\,) & \\
&  V_\zeta \ot V_\zeta & \quad \longrightarrow \quad & \K \\
\endmatrix $$

\n so that 

$$(x\cdot y |z) = (x|\jmath(y \cdot \jmath(z)))
= (\jmath(x) | y \cdot \jmath(z))  \tag
4.8$$ 

\n for all $x,y,z$ in the quantum octonions.   At $q = 1$
this reduces to the identity 

$$(x\cdot y |z)_{{}_1} = (x|z \cdot \o y)_{{}_1}. \tag 4.9 $$ 

\n satisfied by the octonions.
\endproclaim
\m
{\bf Proof.} Recall we have shown that  $\jmath \circ p \circ (\jmath \ot \jmath) :
V_\eta \ot V_\zeta \rightarrow V$ is a $\uqd$-module homomorphism, which
is equivalent to saying

$$a \circ \jmath \circ p \circ (\jmath \ot \jmath) = \jmath \circ p \circ (\jmath \ot
\jmath)\Big((\eta \ot \zeta )\Delta(a)\Big)$$

\n for all $a \in \uqd$.  This implies

$$\aligned
\zeta(a) \circ \jmath \circ p \circ (\jmath \ot \jmath) & = \jmath \circ p \circ
(\jmath
\ot
\jmath)\Big((\eta \ot \zeta )\Delta(\zeta(a))\Big)\\
& = \jmath \circ p \circ
(\jmath
\ot
\jmath)\Big((\eta \ot \zeta )(\zeta \ot \zeta)\Delta(a)\Big)\\
& = \jmath \circ p \circ
(\jmath
\ot
\jmath)\Big((\eta \zeta \ot \text{id} )\Delta(a)\Big)\\ 
\endaligned$$

\n so that $\jmath \circ p \circ
(\jmath
\ot
\jmath): V_{\eta \zeta} \ot V \rightarrow V_\zeta$ is
a $\uqd$-module homomorphism.   From this we deduce that

$$\jmath \circ p \circ
(\text{id} \ot \jmath) = \jmath \circ p \circ
(\jmath
\ot
\jmath)\circ (\jmath \ot \text{id}): V_\eta \ot V \rightarrow V_\zeta$$

\n is a $\uqd$-module homomorphism as well.   Thus, the mappings in the
above diagram are scalar multiples of one another.  The scalar can be
discovered by computing the values of them on $x_1 \ot x_{-3} \ot x_{-2}$:

$$(x_1 \cdot x_{-3} | x_{-2}) = q^{-\frac{3}{2}}(x_2|x_{-2}) = 
\frac {-q^{-\frac{3}{2}}q^{-2}}{ q^3 + q^{-3} } = \frac {-q^{-\frac{7}{2}}}{ q^3 +
q^{-3} }
$$

$$(x_1 | \jmath\big (x_{-3}\cdot\jmath(x_{-2})\big)) = (x_1 |
\jmath(q^{-\frac{1}{2}}x_{-1})) = \frac {-q^{-\frac{1}{2}}q^{-3}}{ q^3 + q^{-3} }= \frac
{-q^{-\frac{7}{2}}}{q^3 + q^{-3} }.$$

\n Hence, they agree and we have the result. The fact that $\jmath$ is
an isometry gives last equality
in (4.8).  \qed 
\b
\proclaim{Proposition 4.10} $(x \cdot y | z) = (\jmath(y) | \jmath(z)
\cdot K_{2\rho}^{-1} x)$ for all $x,y,z \in \Bbb O_q = (V,\cdot)$, which at $q = 1$ 
gives

$$(x \cdot y | z)_{{}_1} = (y| \overline x \cdot z)_{{}_1} \tag 4.11$$

\n for the octonions. 
\endproclaim
\b
{\bf Proof.}  This is an immediate consequence of the calculation,

$$\aligned 
(x \cdot y | z) {\eqa}  
 (x | \jmath(y \cdot \jmath(z))) & {\eqb} (\jmath(y 
\cdot \jmath(z)) | K_{2\rho}^{-1} x)\\
& = (y \cdot\jmath(z) | \jmath(K_{2\rho}^{-1} x)) \eqc (\jmath(y)
| \jmath(z)\cdot K_{2\rho}^{-1} x), \\
\endaligned $$

\n where (1) is from Proposition 4.7, (2) uses (2.20), and (3) is Proposition 4.7. 
 \qed
\b
Note $K_{2\rho}^{-1}$ acts on $\Bbb O_q$ as an algebra
automorphism because it is group-like.  
\b
\proclaim{Proposition 4.12}  The following diagram commutes

$$\matrix &  & \text{id}\ot (\,|\,) \ot \text{id} &   \\
& V_\eta \ot {\Cal S}^2(V_\zeta) \ot V_\eta  & \quad \longrightarrow \quad & V_\eta \ot
V_\eta  \\ & & &  \\
\Big(\jmath \circ p \circ (\jmath \ot \jmath)\Big)\ot p&  \downarrow 
& \quad & \downarrow
(\,|\,)\\ & & (\,|\,) & \\
& V \ot V  & \quad \longrightarrow \quad & \K. \\
\endmatrix $$ 
 
\n When  $q = 1$, this gives the identity  
$$(x\cdot y |x \cdot z)_{{}_1} = (x|x)_{{}_1}(y|z)_{{}_1} \tag 4.13$$

\n for the octonions.  
\endproclaim
\b
{\bf Proof.}  First, it is helpful to observe that

$$\aligned
\Hom_{\uqd}(V_\eta \ot {\Cal S}^2(V_\zeta) \ot V_\eta,\K) & \cong
\Hom_{\uqd}({\Cal S}^2(V_\zeta),V_\eta \ot V_\eta) \\
& \hskip 1 truein = \Hom_{\uqd}(L(0),L(0)),\\
\endaligned
$$

\n so these spaces have dimension 1.  We evaluate the mappings on
$x = x_{-1} \ot (x_1 \ot x_{-1} + x_{-1} \ot x_1)\ot x_1$:

$$\aligned 
(\,|\,) \circ \Big(\text{id}\ot (\,|\,) \ot \text{id}\Big)(x) & = (x_{-1}|x_1) = \frac
{q^{3}}{q^{3} + q^{-3}} \\
(\,|\,) \circ \Big(\big(\jmath \circ p \circ (\jmath \ot \jmath)\big)\ot p\Big)(x) 
& = (\jmath(x_{-1}\cdot x_1) | x_{-1}\cdot x_1)  \\
&  = -q^{3}(x_{-4}|x_4) 
= \frac {q^{3}}{q^{3}
+ q^{-3}}. \qed  \\ 
\endaligned $$  
\b
\proclaim{Proposition 4.14}  The following diagram commutes

$$\matrix &  & \text{id}\ot (\,|\,) \ot \text{id} &   \\
& V_\zeta \ot {\Cal S}^2(V_\eta) \ot V_\zeta  & \quad \longrightarrow \quad & V_\zeta \ot
V_\zeta 
\\ & & &  \\
p \ot \Big(\jmath  \circ p \circ (\jmath \ot \jmath)\Big) & \downarrow 
& \quad & \downarrow
(\,|\,)\\ & & (\,|\,) & \\
& V \ot V  & \quad \longrightarrow \quad & \K \\
\endmatrix $$

\n  In particular, at $q = 1$
this gives the following identity satisfied by the octonions:

$$(y\cdot x |z \cdot x)_{{}_1} = (y|z)_{{}_1}(x|x)_{{}_1} \tag 4.15$$ 
\endproclaim
\b
{\bf Proof.}  The proof is virtually identical to that of Proposition 4.12.
We check the maps on  $x = x_{-1} \ot (x_1 \ot x_{-1} + x_{-1} \ot x_1)\ot x_1$:
$$\aligned 
(\,|\,) \circ \Big(p \ot \big(\jmath \circ p \circ (\jmath \ot \jmath)\big)\Big) (x) 
& = (x_{-1}\cdot x_1|\jmath(x_{-1}\cdot x_1) )  
= (\jmath(x_{-1}\cdot x_1)|x_{-1}\cdot x_1) \\
& = -q^{3}(x_{-4}|x_4) 
= \frac {q^{3}}{q^{3}
+ q^{-3}},  \\ 
\endaligned $$

\n which equals  $(\,|\,) \circ \big(\text{id}\ot (\,|\,) \ot \text{id}\big)(x)$.  \qed
\b
\subhead {$r$-algebras}\endsubhead
\b
The notion of an $r$-algebra (that is, an algebra equipped with a Yang-Baxter operator)
arises in Manin's work [M] on noncommutative geometry.
Noteworthy examples of $r$-algebras are 
Weyl and Clifford algebras, noncommutative tori, certain
universal enveloping algebras, and quantum groups (see for
example, Baez [B]).  
\b
\proclaim {Definition 4.16}  Suppose $A$ is an algebra
over a field $\F$  with multiplication $p: A \ot
A \rightarrow A$.  Assume 
$R \in End(A^{\ot 2})$ is a Yang-Baxter operator, i.e. 
$R_{1,2}R_{2,3}R_{1,2} =
R_{2,3}R_{1,2}R_{2,3}$, where $R_{1,2} = R \ot \text{id} \in \End(A^{\ot 3})$
and $R_{2,3} = \text{id} \ot R \in  \End(A^{\ot 3})$.   Then $A$ is
said to be  an $r$-algebra if

$$\aligned
& R \circ (p \ot \text {id}_A) = (\text {id}_A \ot p)\circ R_{1,2}\circ R_{2,3} \\
& R \circ (\text {id}_A \ot p) = (p \ot \text {id}_A) \circ R_{2,3} \circ R_{1,2}. \\
\endaligned \tag 4.17
$$

\endproclaim

When $A$ has a unit element 1, the Yang-Baxter operator $R$ is also required
to satisfy  $R(1 \ot x) = x \ot 1, \quad  R(x \ot 1) = 1 \ot x \quad \text
{for all}\; x \in A.$ 
\b
Our quantum octonions satisfy $r$-algebra properties similar to those in (4.17) with
respect to the R-matrix of $\uqd$.
\b
\proclaim{Proposition 4.18} (a) In $\Hom_{\uqd}(V_\zeta \ot V_\eta \ot V, V \ot V)$ we have
$$\check R_{V,V} \circ (p \ot \text {id}) = (\text {id} \ot p)\circ 
\check R_{V_\zeta,V}\circ \check R_{V_\eta,V}.$$
\s
\item{(b)}{} In $\Hom_{\uqd}(V \ot V_\zeta \ot V_\eta, V \ot V)$ we have

$$\check R_{V,V} \circ (\text {id} \ot p) = (p \ot \text {id}) \circ 
\check R_{V,V_\eta} \circ \check R_{V,V_\zeta}.$$ \endproclaim
\b
{\bf Proof.} Suppose that $\pi_1: V \ot V \rightarrow L(2\omega_1)$,
$\pi_2: V \ot V \rightarrow L(\omega_2),$
and  $\pi_3: V \ot V \rightarrow L(0)$ are the projections of
$V \ot V$ onto its irreducible summands. Let $T_\ell$
and $T_r$ denote the transformations on the left and right
in (a).  We compute $T_\ell$ and $T_r$ on several elements of
$V_\zeta \ot V_\eta \ot V$ which have been chosen so that
$\Theta$ acts as the identity.  We use the fact that
$\eta(\e_1) = 1/2(\e_1+ \e_2 + \e_3 + \e_4)$,
$\zeta(\e_1) = 1/2(\e_1+ \e_2 + \e_3 - \e_4)$,
$\eta(\e_2) = 1/2(\e_1+ \e_2 - \e_3 - \e_4)$, and
$\zeta(\e_3) = 1/2(\e_1- \e_2 + \e_3 + \e_4)$.

We begin by evaluating the action of $T_\ell$ and $T_r$  on
$x_{-1} \ot x_1 \ot x_{-1}$:

$$\aligned\Big(\check R_{V,V} \circ (p \ot \text {id})\Big)(x_{-1} \ot x_1 \ot x_{-1})
& = -q^{\frac{3}{2}}\check R_{V,V}(x_{4} \ot x_{-1}) = 
-q^{\frac{3}{2}}(x_{-1} \ot x_4) \\
\Big( (\text {id} \ot p) \circ 
\check R_{V_\zeta,V}\circ \check R_{V_\eta,V}\Big)(x_{-1} \ot x_1 \ot x_{-1})
& = q^{\frac 12}\Big((\text {id} \ot p) \circ 
\check R_{V_\zeta,V}\Big)(x_{-1} \ot x_{-1} \ot x_1) \\
& = 
(\text {id} \ot p)(x_{-1} \ot  x_{-1} \ot x_1)\\
& = -q^{\frac{3}{2}}(x_{-1} \ot x_4).\\
\endaligned$$

\n Now $(x_{-1}|x_4) = 0$ implies that $x_{-1} \ot x_4$
lies in $L(2\omega_1) \oplus L(\omega_2)$. The vector
$x_{-1} \ot x_4$ is not an eigenvector of $\check R_{V,V}$ 
because  $\check R_{V,V}(x_4 \ot x_{-1}) = x_{-1} \ot x_4$.
Therefore, $\pi_1(x_{-1} \ot x_4) \neq 0$ and $\pi_2(x_{-1} \ot x_4)
\neq 0$.   Since 
$\dim \Hom_{\uqd}(V_\zeta \ot V_\eta \ot V, L(2\omega_1)) = 1
= \dim \Hom_{\uqd}(V_\zeta \ot V_\eta \ot V, L(\omega_2)$,
we see from the above calculation that 
$\pi_1 T_\ell = \pi_1 T_r$  
and $\pi_2 T_\ell = \pi_2 T_r$.
\m
Consider the action of
$T_\ell$ and $T_r$ on $x_3 \ot x_2 \ot x_{-1}$:  

$$\aligned\Big(\check R_{V,V} \circ (p \ot \text {id})\Big)(x_3 \ot x_2 \ot x_{-1})
& = q^{\frac{1}{2}}\check R_{V,V}(x_{1} \ot x_{-1}) = 
q^{\frac{3}{2}}(x_{-1} \ot x_1) \\
\Big( (\text {id}\ot p) \circ 
\check R_{V_\zeta,V}\circ \check R_{V_\eta,V}\Big)(x_3 \ot x_2 \ot x_{-1})
& = q^{\frac 12}\Big((\text {id}\ot p) \circ 
\check R_{V_\zeta,V}\Big)(x_3 \ot x_{-1} \ot x_2) \\
& = 
q(\text {id}\ot p)(x_{-1} \ot x_3 \ot x_2)\\
& = q^{\frac{3}{2}}(x_{-1} \ot x_1).\\
\endaligned$$   

\n The projection $\pi_3$ comes from the form 
$(\,|\,)$, and so it is nonzero when
it acts on $x_{-1} \ot x_1$.  Consequently, applying $\pi_3$ to
both sides of 
$T_\ell(x_3 \ot x_2 \ot x_{-1}) = T_r(x_3 \ot x_2 \ot x_{-1})$ 
and using $\dim \Hom_{\uqd}(V_\zeta \ot V_\eta \ot V, L(0)) = 1$,
we obtain $\pi_3 T_\ell = \pi_3 T_r$.  Therefore,
$T_\ell = (\pi_1 + \pi_2 + \pi_3)T_\ell = (\pi_1 + \pi_2 + \pi_3)T_r = T_r.$
The calculations for part (b) are almost identical and are left
as an exercise for the reader.  \qed    
\b
\s

\head {\S 5. Properties of the Quantum Para-Hurwitz Algebra} \endhead 
\b
Recall that the quantum para-Hurwitz algebra $\Bbb P_q = (V,\ast)$ is obtained
from identifying the underlying vector spaces of $V_\theta$, 
$V_{\theta^2}$
and $V$, where $V$ is the natural representation of $\uqd$ and $\theta$
is the graph automorphism $\theta = (1\,4\,3)$, and using the
multiplication  $p^*: V_\theta \ot V_{\theta^2} \rightarrow V$, $x \ot y \mapsto x \ast y
= \jmath(x)\cdot \jmath(y)$.  The resulting product is displayed in Table 2.7.   
As we saw in the last chapter, the quantum octonion algebra enjoys 
very nice composition
properties relative to the bilinear $(\,|\,)$.   It is
known that the para-Hurwitz algebra $\Bbb P$ (that is
the $q = 1$ version of $\Bbb P_q$), which is obtained from
the octonions using the product $x \ast y = \overline x \cdot \overline y$,
has an associative bilinear form admitting composition.  
Here we establish similar results for the
bilinear form on $\Bbb P_q = (V,\ast)$ 
\b
\proclaim{Proposition 5.1} For the quantum para-Hurwitz algebra $\Bbb P_q = (V,\ast)$,

$$(x \ast y| z) = (x | y \ast z) \quad \quad \text {and} \quad \quad (x
\ast y | z) = (y | z \ast K_{2\rho}^{-1} x)$$

\n for all $x,y,z \in \Bbb P_q$.   At $q = 1$ these relations reduce to the
well-known properties of the para-Hurwitz algebra $(x \ast y| z)_{{}_1} = (x | y \ast
z)_{{}_1}$ and  $(x
\ast y | z)_{{}_1} = (y | z \ast x)_{{}_1}$.
\endproclaim
\b
{\bf Proof.}  The proof amounts to using the fact that $\jmath$ is an isometry together
with the fact that $\jmath$ commutes with $K_{2\rho}^{-1}$.  So we have

$$(x \ast y|z) = (\jmath(x)\cdot \jmath(y) | z) = 
\cases (x | \jmath(y)\cdot \jmath(z)) = (x | y \ast z) \\
(y | \jmath(z) \cdot K_{2\rho}^{-1} \jmath(x)) =
(y | \jmath(z) \cdot \jmath (K_{2\rho}^{-1} x)) = 
(y | z \ast K_{2\rho}^{-1} x). \\
\endcases \qed $$
\b
\proclaim {Proposition 5.2} Ignoring the $\uqd$-module structures we have

$$\matrix &  & \text{id}\ot (\,|\,) \ot \text{id} &   \\
& V \ot {\Cal S}^2(V) \ot V  & \quad \longrightarrow \quad & V \ot V
\\ & & &  \\
(\jmath \circ p^* \circ (\jmath\ot \jmath)) \ot p^*  & \downarrow 
& \quad & \downarrow
(\,|\,)\\ & & (\,|\,) & \\
& V \ot V  & \quad \longrightarrow \quad & \K \\
\endmatrix $$
\s
$$\matrix &  & \text{id}\ot (\,|\,) \ot \text{id} &   \\
& V \ot {\Cal S}^2(V) \ot V  & \quad \longrightarrow \quad & V \ot V
\\ & & &  \\
p^* \ot \Big(\jmath \circ p^* \circ (\jmath \ot \jmath) \Big)  & \downarrow 
& \quad & \downarrow
(\,|\,)\\ & & (\,|\,) & \\
& V \ot V  & \quad \longrightarrow \quad & \K \\
\endmatrix$$ 
\endproclaim
\b
{\bf Proof.}  These diagrams commute when $p$ is used in place of $p^*$.  
Now when  $p^* = p \circ (\jmath \ot \jmath)$ is taken, we have

$$\aligned (\,|\,) \circ (\jmath \circ p^* \circ (\jmath\ot \jmath)) \ot p^*
& = (\,|\,) \circ \Big ((\jmath \circ p) \ot (p \circ (\jmath \ot \jmath))\Big) \\
& \eqa (\,|\,) \circ p \ot (j \circ p \circ (\jmath\ot \jmath)) \\
& = (\,|\,) \circ \text{id}\ot (\,|\,) \ot \text{id} \\ \endaligned $$

\n where (1) follows from the fact that $\jmath$ is an isometry.  The
other diagram can be shown to commute in a similar way.  \qed  
\b 
\proclaim {Proposition 5.3 } The following diagram commutes and
shows  at $q = 1$ that the relation $(x\cdot y)\cdot x =
(x|x)_{{}_1}y = x \cdot(y \cdot x)$ holds for the para-Hurwitz
algebra.

$$
\thickmuskip2mu
\medmuskip2mu
\matrix V_{\theta^2} \ot V & \lmap & {\Cal S}^2(V_{\theta^2})\ot V_\theta = {\Cal S}^2(V)
\ot V = {\Cal S}^2(V_\theta)\ot V_{\theta^2} & \ppmap & V_\theta \ot V   \\
& & &   &     \\
 \downarrow p^* & &(\,|\,)\ot \text{id}\downarrow & 
\hskip-10pt
\jmath \circ p^* \circ
(\jmath
\ot
\jmath)\hskip-10pt
& \downarrow \\
 & & &   &    \\
V_\theta &  = &   V &   = &   V_{\theta^2}  \\
\endmatrix $$ 
\endproclaim 
\b
{\bf Proof.} Again it suffices to check the maps on the element $x =
(x_1 \ot x_{-1} + x_{-1} \ot x_1) \ot x_1$, and for it we have

$$
\aligned \Big (\jmath \circ p^* \circ (\jmath \ot \jmath)\Big)
\circ (\text{id} \ot p^*)(x) & =  \jmath \circ p^* \circ (\jmath \ot \jmath)
(x_1 \ot (-q^{\frac{3}{2}})x_4)  \\
& = -\jmath(x_1) = x_1 = (\,|\,) \ot \text{id}(x), \quad \text{while}\\
\endaligned $$

$$
p^* \circ \Big (\text{id} \ot \big(\jmath  \circ p^* \circ (\jmath \ot
\jmath)\big)\Big)(x) = p^*(x_1 \ot q^{\frac{3}{2}}x_{-4}) = x_1 = (\,|\,) \ot \text{id}(x).
\qed $$

\b 
\b
 
\head {\S 6. Idempotents and Derivations of Quantum Octonions} \endhead
\b
It is apparent from Table 3.9 that $e_1 = q^{-\frac{1}{2}}x_4$ and
$e_2 = -q^{\frac {1}{2}}x_{-4}$ are orthogonal idempotents in the quantum octonion
algebra $\Bbb O_q = (V,\cdot)$.   Their sum $e = e_1 + e_2$ acts as a left and right
identity element on the span of $e_1,e_2, x_j,x_{-j}, j = 2,3$.  Moreover, the
relations 
$e\cdot x_1 = q^{2}x_1$, $x_1\cdot e = q^{-2}x_1$, $e \cdot x_{-1} = q^{-2}x_{-1}$,
and $x_{-1} \cdot e = q^{2} x_{-1}$ hold in $\Bbb O_q$.   Hence, at $q = 1$ or $-1$, the
element $e$ is an identity element.  Relative to the basis
$e_1,e_2, x_{\pm j}, j = 1,2,3$, the left and right multiplication operators
$L_{e_1}, R_{e_1},L_{e_2}, R_{e_2}$ can be simultaneously diagonalized. 
The operators $L_{e_1}, R_{e_1}$ determine a Peirce decomposition
of $\Bbb O_q$,  

$$\Bbb O_q = \oplus_{(\ga,\delta)} (\Bbb O_q)_{\ga,\delta} \quad \text {where} \quad
(\Bbb O_q)_{\ga,\delta} = \{ v
\in \Bbb O_q \mid e_1 \cdot v = \ga v, \; \text {and} \; v\cdot e_1 = \delta v\}. \tag
6.1$$
\m
\n 
Now $(\Bbb O_q)_{0,q^{-2}} = \K x_1$, $(\Bbb O_q)_{1,0} = \K x_2 + \K x_3$, 
$(\Bbb O_q)_{1,1} = \K e_1$, $(\Bbb O_q)_{0,0} = \K e_2$, $(\Bbb O_q)_{0,1} =
\K x_{-2} +
\K x_{-3}$, and $(\Bbb O_q)_{q^{-2},0} = \K x_{-1}$.   The idempotent $e_2$ gives a
similar Peirce decomposition.
\b
Suppose $d$ is a derivation of $\Bbb O_q$ such that $d(e_1) = 0 = d(e_2)$.
It is easy to see that $d$ must map each Peirce space $(\Bbb O_q)_{\ga,\delta}$ into
itself.  Thus, we may suppose that
\m 
$$\gathered
d(x_1) = b_1 x_1 \quad \quad d(x_{-1}) = b_{-1} x_{-1}  \\
d(x_2) = c_{2,2}x_2 + c_{3,2} x_3 \quad \quad 
d(x_3) = c_{2,3}x_2 + c_{3,3} x_3 \\
d(x_{-2}) = d_{2,2}x_{-2} + d_{3,2} x_{-3} \quad
\quad  d(x_{-3}) = d_{2,3}x_{-2} + d_{3,3} x_{-3}. \\  
\endgathered  \tag 6.2
$$

\n Applying $d$ to the relation $x_1 \cdot x_{-1} = q^{-\frac {3}{2}} x_{-4} =
-q^{-2}e_2$ (or to the relation $x_{-1} \cdot x_1 = -q^{\frac {3}{2}}x_4
= -q^{2} e_1$)  shows that $b_{-1}
= -b_1$.   Now from
$x_k \cdot x_k = 0$, we deduce that $c_{j,k} = 0$ for $j \neq k$, and similarly 
from $x_{-k} \cdot x_{-k} = 0$ it follows that
$d_{j,k} = 0$.  There are relations $x_1 \cdot x_{-j} = \xi x_k$, $x_{-j} \cdot x_1 =
\xi x_k$, $x_{-1} \cdot x_j =
\xi x_{-k}$, and $x_j \cdot x_{-1} = \xi x_{-k}$ (where $j \neq k$ and $\xi$ in each one
indicates some appropriate scalar that can be found in
Table 3.9), and applying
$d$ to those equations gives $b_1 + d_{j,j} = c_{k,k}$ and $b_{-1} + c_{j,j} = d_{k,k}$.
The second is equivalent to the first since $b_{-1} = -b_1$.  
From  $x_{\pm j} \cdot x_{\pm k} = \xi x_{\pm 1}$ we see that
$c_{2,2} + c_{3,3} = b_1$ and $d_{2,2} + d_{3,3} = -b_1$.  Finally, applying
$d$ to $x_j \cdot x_{-j} = \xi e_{1}$ or
$x_{-j} \cdot x_j = \xi e_2$, we obtain $d_{j,j} = -c_{j,j}$
for $j = 2,3$. Every nonzero product between basis elements in the quantum octonions
is one of these types, or it is one involving the idempotents $e_1,e_2$, so
we have determined all the possible relations.   Consequently, we can conclude that
for any derivation $d$ such that
$d(e_1) = 0 = d(e_2)$, there are scalars $b = b_1$ and $c = c_{2,2}$ so that the
following hold:

$$\gathered 
d(e_1) = 0 = d(e_2), \\
 d(x_1) = bx_1 \quad \quad d(x_{-1}) = -b x_{-1} \\
d(x_2) = c x_2 \quad \quad d(x_{-2}) = -c x_{-2} \\
d(x_3) = (b-c)x_3 \quad \quad  d(x_{-3}) = -(b-c) x_{-3}. \\
\endgathered \tag 6.3
$$ 

\n The same calculations show that any transformation defined by (6.3) where
$b,c$ are arbitrary scalars is a derivation of $\Bbb O_q$.  

Now suppose that $d$ is an arbitrary derivation of $\Bbb O_q$, and  
let $d(e_1) = a_0 e_1 + a_0' e_2 + \sum_{j = -3,j \neq 0}^3 a_j x_j$ and
$d(e_2) = b_0 e_1 + b_0' e_2 + \sum_{j = -3,j \neq 0}^3 b_j x_j$.  Computing
$d(e_1) \cdot e_1 + e_1 \cdot d(e_1) = d(e_1)$, we obtain

$$2a_0 e_1 + q^{-2}a_1 x_1 +  a_2 x_2 + a_3 x_3 + q^{-2}a_{-1}x_{-1} + 
a_{-2} x_{-2} + a_{-3} x_{-3} = d(e_1),$$

\n which shows that $a_0 = a_0' = a_{-1} = a_1 = 0$.   A similar
calculation with the relation $d(e_2^2) = d(e_2)$ gives $b_0 = b_0' = b_{-1} = b_1 = 0$.
From $d(e_1)\cdot e_2 + e_1 \cdot d(e_2) = 0$ it follows that $b_2 = -a_2$ and $b_{3} = -
a_3$, and from $d(e_2) \cdot e_1 + e_2 \cdot d(e_1) = 0$, the relations $b_{-2} =
-a_{-2}$,
$b_{-3} = -a_{-3}$ can be determined.  

Now $e_1 \cdot d(x_1) = - d(e_1)\cdot x_1$ shows that the coefficients of $x_2$ and $x_3$
in $d(x_1)$ are $q^{\frac {3}{2}}a_{-3}$ and  $q^{\frac {3}{2}}a_{-2}$ respectively;
while $d(x_1) \cdot e_2 = -x_1 \cdot d(e_2)$ gives that those coefficients are
$-q^{-\frac {3}{2}}b_{-3} = q^{-\frac {3}{2}}a_{-3}$  and
$-q^{-\frac {3}{2}}b_{-2} = q^{-\frac {3}{2}}a_{-2}$.  Thus, when
$q$ is not a cube root of unity, we obtain that $a_{-2} = 0 = a_{-3}$.
Likewise  the relations $a_2 = 0 = a_3$
can be found from analyzing $d(x_{-1})\cdot e_1 = -x_{-1} d(e_1)$ and 
$e_2 \cdot d(x_{-1}) = - d(e_2)\cdot x_{-1}$.  Consequently, when $q$ is not  
a cube root of unity,  every derivation of $\Bbb O_q$ must annihilate
$e_1$ and $e_2$, and so must be as in (6.3) for suitable scalars $b,c$.  To summarize we
have
\b
\proclaim {Theorem 6.4} The derivation algebra of the quantum octonion algebra $\Bbb O_q =
(V,\cdot)$ is two-dimensional.  Every derivation $d$ of $\Bbb O_q$ is given by  (6.3) for some
scalars $b,c \in \K$.  
\endproclaim
\b
\subhead {Simplicity of the quantum octonions}\endsubhead
\b
The idempotents also provide a proof of the following
\b
\proclaim {Theorem 6.5} The quantum octonion algebra $\Bbb O_q = (V,\cdot)$ is simple. 
\endproclaim
\b
{\bf Proof.}  Let $I$ be a nonzero ideal of $\Bbb O_q$.  Since $e_1 I \subseteq I$
and $I e_1 \subseteq I$,  it follows that

$$I = \oplus_{\ga,\delta} I_{\ga,\delta},$$

\n where $I_{\ga,\delta} = \{x \in I \mid e_1 \cdot x = \ga x,$ and $x \cdot e_1
= \delta x\}$.  Assume first $x_i \in I$ for some $i = \{\pm 1, \pm 2, \pm 3, \pm
4\}$.  
If we multiply $x_i$ by $x_{-i}$ when $i = \pm 1, \pm 2,$ or $\pm 3$, we see
that $e_1$ or $e_2$ belongs to $I$.  Therefore, either  $e_1 \Bbb O_q + \Bbb O_q e_1 
\subseteq I$ or $e_2 \Bbb O_q + \Bbb O_q e_2 \subseteq I$, and as a consequence,
$I = \Bbb O_q$.   It remains to consider the possibility
that some nonzero linear combination $a x_2 + b x_3$ or $a x_{-2} + b x_{-3}$
is in $I$.  Multiplying by one of the elements $x_2, x_3, x_{-2}$, or $x_{-3},$
we obtain that $x_1$ or $x_{-1}$ is in $I$, and that implies $I = \Bbb O_q$
as before.  Consequently, $\Bbb O_q$ is simple.  \qed  
\b   
\b
\head {\S 7. Quantum Quaternions} \endhead
\b
The subalgebra $\Bbb H_q$ of $\Bbb O_q$ generated by $x = x_2$ and $y = x_{-2}$ is
4-dimensional and has $x,y,e_1,e_2$ as a basis, where $e_1 = q^{-\frac{1}{2}}x_4$ and
$e_2 = -q^{\frac {1}{2}}x_{-4}$ are the orthogonal idempotents of
Section 6.   The sum $e = e_1 + e_2$ is the identity element of $\Bbb H_q$, and
multiplication in this algebra is given by the following: 

$$\table  {\bf \it Table 7.1}
\hfill ! \hfill \hfill ! \hfill \hfill ! \hfill \hfill ! \hfill  \r 
  \hfill  $\cdot$ \hfill ! \hfill $e_1$ \hfill ! \hfill $e_2$ \hfill ! \hfill $x$ \hfill
! \hfill $y$ \hfill \rr 
\hfill  $e_1$ \hfill ! \hfill $e_1$ \hfill ! \hfill 0 \hfill ! \hfill $x$ \hfill ! \hfill
$0$ \hfill \rr
\hfill  $e_2$ \hfill ! \hfill 0 \hfill ! \hfill $e_2$ \hfill ! \hfill 0 \hfill ! \hfill
$y$ \hfill  \rr
\hfill  $x$ \hfill ! \hfill 0 \hfill ! \hfill $x$ \hfill ! \hfill 0 \hfill ! \hfill
$q^{-1}e_1$ \hfill   \rr
\hfill  $y$ \hfill ! \hfill $y$ \hfill ! \hfill 0 \hfill ! \hfill $q e_2$ 
\hfill ! \hfill
$0$ \hfill  \caption{}
$$

\n It is easy to see that at $q = 1$ the algebra above is isomorphic
to the algebra $M_2(\K)$  of  $2 \times 2$ matrices over $\K$ via the
mapping $e_1 \mapsto E_{1,1}$, $e_2 \mapsto E_{2,2}$, $x \mapsto E_{1,2}$,
$y \mapsto E_{2,1}$, which sends basis elements to standard matrix units.
Since the split quaternion algebra is isomorphic to $M_2(\K)$,
this algebra reduces to the quaternions at $q = 1$.  Another way to see
this is to observe that the bilinear form on $\Bbb O_q$, when restricted
to $\Bbb H_q$, satisfies (3.6) at $q = 1$. Thus, at $q = 1$ we have 
a 4-dimensional unital algebra with a nondegenerate bilinear form admitting composition,
and so it must be the quaternions. 
\m
The algebra $\Bbb H_q$ is neither associative nor flexible if $q \neq \pm 1$  since

$$(x \cdot y) \cdot x = q^{-1} x \neq q x = x \cdot (y \cdot x).$$

\n However, it is {\it Lie admissible}.  Indeed, if
$h = q^{-1}e_1-q e_2$, then

$$[h,x] = (q+q^{-1})x, \quad  [h,y] = -(q+q^{-1})y, \quad \text{and} \quad [x,y] = h.$$
  
\n Therefore, the elements 
$x, \, \displaystyle{\frac {2}{q+q^{-1}}}y, \, \displaystyle{\frac {2}{q+q^{-1}}} h$  
determine 
a standard basis for $sl_2$.  Since $e = e_1 + e_2$ is the identity, we see that
$\Bbb H_q$ is the Lie algebra $gl_2$ 
under the commutator product $[a,b] = ab - ba$.  
\m
The subalgebra generated by $x_3$ and $x_{-3}$ is isomorphic to the
algebra $\Bbb H_q$.  However, that is not true if we consider the
subalgebra $\Bbb H'_q$ generated by $x_1$ and $x_{-1}$.  Let $u = x_{-1}$
and $v = x_1$.  Then multiplication in  $\Bbb H'_q$ is given by  
\m
$$\table  {\bf \it Table 7.2}
\hfill ! \hfill \hfill ! \hfill \hfill ! \hfill \hfill ! \hfill  \r 
  \hfill  $\cdot$ \hfill ! \hfill $e_1$ \hfill ! \hfill $e_2$ \hfill ! \hfill $u$ \hfill
! \hfill $v$ \hfill \rr 
\hfill  $e_1$ \hfill ! \hfill $e_1$ \hfill ! \hfill 0 \hfill ! \hfill $q^{-2}u$ \hfill !
\hfill
$0$ \hfill \rr
\hfill  $e_2$ \hfill ! \hfill 0 \hfill ! \hfill $e_2$ \hfill ! \hfill 0 \hfill ! \hfill
$q^2v$ \hfill  \rr
\hfill  $u$ \hfill ! \hfill 0 \hfill ! \hfill $q^2 u$ \hfill ! \hfill 0 \hfill ! \hfill
$q^2e_1$ \hfill   \rr
\hfill  $v$ \hfill ! \hfill $q^{-2}v$ \hfill ! \hfill 0 \hfill ! \hfill $q^{-2}e_2$ 
\hfill ! \hfill
$0$ \hfill  \caption{}
$$

There is no identity element in $\Bbb H'_q$.  But at  $q = 1$, the element $e_1 + e_2$ 
becomes an identity element, and the algebra is isomorphic to
$M_2(\K)$.  The bilinear form at $q =1$ has the
composition property, so the algebra $\Bbb H_q'$  too might rightfully be termed a quantum
quaternion algebra.  
The elements $h' = q^2 e_1 -q^{-2} e_2, u, v$ determine a standard basis
for $sl_2$:

$$[h',u] = 2u, \quad \quad [h',v] = -2v, \quad \quad [u,v] = h'.$$

\n Since $[q^2 e_1 +q^{-2} e_2,h'] = [q^2 e_1 +q^{-2} e_2,u] = 
[q^2 e_1+q^{-2} e_2,v] = 0$, the algebra  $\Bbb H_q'$ under the
commutator product is isomorphic to $gl_2$.  It is neither 
associative nor flexible because

$$\big((e_1+v)\cdot e_1\big)\cdot (e_1+v) = e_1 + q^{-4} v \neq 
e_1 + q^{-2} v = (e_1+v)\cdot \big(e_1\cdot (e_1+v)\big).$$
 
\b
\head {\S 8. Connections with $U_q$(B$_3$) and $U_q$(G$_2$)} \endhead
\b
In this section we assume that the field $\K$ has characteristic zero. We consider the
subalgebras of
$\uqd$ that are described as the fixed points of graph automorphisms by

$$\aligned
U_{\zeta,\eta} & = \{ a \in \uqd \mid \zeta(a) = a = \eta(a)\} \\
U_\delta & = \{a \in \uqd \mid \delta(a) = a \}, \\
\endaligned \tag 8.1 $$

\n where $\delta = (3\,4) = \zeta \eta \zeta$.  Then $U_\delta \supseteq U_{\zeta,\eta}$.
In the non-quantum case, the enveloping algebra $U($G$_2) \subseteq U_{\zeta,\eta} \subseteq
U($D$_4)$ and $U($B$_3) \subseteq U_\delta \subseteq U($D$_4)$.  
We consider the structure of $V,V_\zeta,V_\eta$ as modules for
the subalgebras given in (8.1).  
\b
Now for any $a \in U_\delta$, we have that

$$\jmath(a x) = \jmath(\delta(a) x) = a \jmath (x)$$

\n for all $x \in V$ so that any $a \in U_\delta$ commutes with $\jmath$. 
In particular, the eigenspaces of $\jmath$ are $U_\delta$-invariant.
But the eigenspaces of $\jmath$ are 1-dimensional (spanned by  
 $x_{4}-x_{-4}$), and 7-dimensional (spanned by $\{x_i \mid i = \pm 1,\pm 2, \pm 3\} \cup
\{x_4 + x_{-4}\})$.  Consider the 7-dimensional space under the
action of $U_{\zeta,\eta}$.  Since $F_{\a_2}, F_{\a_1}+F_{\a_3} + F_{\a_4},
E_{\a_2}, E_{\a_1}+E_{\a_3} + E_{\a_4} \in U_{\zeta,\eta}$, it is easy to
check that the 7-dimensional space is an irreducible $U_{\zeta,\eta}$-module,
hence an irreducible $U_\delta$-module as well.  The action of $U_{\zeta,\eta}$
on $V_\zeta$ and $V_\eta$ is the same as its action on $V$.  As a result we
have
\b
\proclaim {Proposition 8.2} The modules $V,V_\zeta \cong V_-,$ and $V_\eta \cong
V_+$ are isomorphic as modules for the subalgebra $U_{\zeta,\eta}$
in (8.1).  They decompose into the sum of  
irreducible $U_{\zeta,\eta}$-modules of dimensions 1 and 7, which are the eigenspaces of
$\jmath$. 
\endproclaim
\b
For $q = 1$, Proposition 8.2 amounts to saying
that $V, V_\zeta, V_\eta$ are isomorphic $U$(G$_2$)-modules
which decompose into irreducible modules of dimension 1 and 7.   

\b
The elements $F_{\a_1}, F_{\a_2}, F_{\a_3}+F_{\a_4}, 
E_{\a_1}, E_{\a_2}, E_{\a_3}+E_{\a_4}$ belong to $U_\delta$, so if
we regard $V_\zeta$ or $V_\eta$ as a module for $U_\delta$, then
we obtain a diagram  
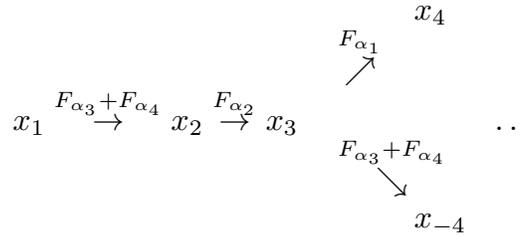
\b
$\hskip 3.4 truein  x_4$

$\hskip 3 truein \fcmap $
 
$\hskip 1.3 truein x_1 \famap x_2 \fbmap x_3 \hskip 1 truein \dots $

$\hskip 3  truein \fdmap$

$\hskip 3.4 truein  x_{-4}$
  
\n which is gotten by superimposing the two diagrams of $V_\zeta$
and $V_\eta$.  Using that, it is easy to argue that both these
modules remain irreducible for $U_\delta$.  Moreover, for any $a \in U_\delta$,
we have

$$\aligned
\jmath(\zeta(a)x) & = \jmath((\zeta \delta)(a)x) = \jmath ((\eta \zeta)(a)x) \\
& = \jmath(\theta^2(a)x) = \eta(a)\jmath(x)\\
\endaligned$$

\n where $\theta = (1 \, 4 \, 3)$.   Consequently, $V_\zeta \cong V_\eta$ as 
$U_\delta$-modules.  
\b
At $q = 1$ these computations show that $V_\zeta$ and $V_\eta$ (the
two spin modules of $U$(D$_4$)) remain irreducible for  $U$(B$_3$), and
in fact they are the spin representation of that algebra. To summarize
we have

\proclaim {Proposition 8.3} When regarded as modules
for the subalgebra $U_\delta$ of $\uqd$,  $V_\zeta \cong V_-,$ and
$V_\eta \cong
V_+$ are irreducible and isomorphic.  \endproclaim
\b
In light of these propositions it is natural to expect that these fixed
point subalgebras of the graph automorphisms might be the quantum groups $U_q(\text{G}_2)$
and $U_q(\text{B}_3)$.  However, that is not true.  In fact the next
``nonembedding'' result  shows that there are no Hopf algebra homomorphisms
of them into $\uqd$ except for the trivial ones given by the counit.  
\b
\proclaim{Proposition 8.4} When $q$ is not a root of unity, then any
Hopf algebra homomorphisms

$$U_q(\text{G}_2) \rightarrow U_q(\text{B}_3), \quad  
U_q(\text{B}_3) \rightarrow U_q(\text{D}_4),  \quad U_q(\text{G}_2) \rightarrow
U_q(\text{D}_4)$$

\n are trivial ($a \mapsto \e(a)1$ for all $a$).  
\endproclaim

\b 
The method we use to establish this result is due originally to Hayashi [H].
We begin with a little background before presenting the proof.  
\b
Suppose $U_q(\g)$ is a quantum group and
$(X,\pi)$, $\pi: U_q(\g) \rightarrow \End(X)$, is a representation of $U_q(\g)$.
This determines another representation $(X, \pi')$, $\pi' =\pi \circ S^2$, of $U_q(\g)$
where $S$ is the antipode.  But it is well-known [Jn, p. 56] that $S^2$
has the following expression:

$$S^2(a) = K_{2\rho}^{-1} a K_{2\rho} \quad \quad \text{for all}\; a \in U_q(\g)$$

\n where $\rho = 1/2 \sum_{\gamma > 0} \gamma = \sum_{\a \in \Pi} c_\a \a$ (the half-sum of
the positive roots of $\g$) and
$K_{2\rho} = \prod_{\a \in \Pi} K_{\a}^{c_\a}$ (where $\Pi$ is the set of simple roots of
$\g$). The two representations $(X,\pi)$ and $(X,\pi')$ are in fact isomorphic by
the following map

$$\phi: (X,\pi) \rightarrow (X, \pi'), \quad \quad  x \mapsto \pi(K_{2\rho}^{-1})x$$

\n because $\phi(\pi(a)x) = \pi(K_{2\rho}^{-1})\pi(a)x = \pi(K_{2\rho}^{-1}a)x
= \pi(S^2(a)K_{2\rho}^{-1})x =
\pi(S^2(a))\pi(K_{2\rho}^{-1})x = \pi'(a)\phi(x)$ for all $a \in U_q(\g)$.      
\b
{\bf Examples 8.5.} Let $\g =$ D$_4$ and let $X = V,V_-,V_+$ (the natural 
or spin representations of $\uqd$).  Here $\rho = 3 \a_1 + 5 \a_2 + 3 \a_3 + 3 \a_4$
so $K_{2\rho}^{-1}= K_{\a_1}^{-6} K_{\a_2}^{-10} K_{\a_3}^{-6} K_{\a_4}^{-6}$. 
Relative to the basis $\{x_1,x_2,x_{3},x_4,x_{-4},x_{-3},x_{-2},x_{-1}\}$
we have as before  

$$\phi_{{\text D}_4} = \phi^-_{{\text D}_4} = \phi^+_{{\text D}_4}
= \text{diag}\{q^{-6},q^{-4},q^{-2},1,1,q^2,q^4,q^6\}.\tag 8.6$$
\b
Next let $\g =$ G$_2$ and let $X$ be the 7-dimensional irreducible
$U_q(\text{G}_2)$-representation.  Here $\rho = 5 \a_1 + 3 \a_2$ so
$K_{2\rho}^{-1}= K_{\a_1}^{-10} K_{\a_2}^{-6}$.  With respect to the
following basis (see [Jn], Chap. 5),

$$\{x_{2\a_1+\a_2}, x_{\a_1+\a_2}, x_{\a_1}, x_0, x_{-\a_1},
x_{-(\a_1+\a_2)}, x_{-(2\a_1+\a_2)}\}, \quad \quad \text{we have} \tag 8.7 $$ 

$$\phi_{\text{G}_2} 
= \text{diag}\{q^{-10},q^{-8},q^{-2},1, q^2,q^8,q^{10}\}.\tag 8.8$$
\b
Finally, assume $X$ is the 7-dimensional natural representation or the 
8-dimensional spin representation of $U_q($B$_3)$.    In this
case $2\rho = 5 \a_1 + 8 \a_2 + 9 \a_3$ so that 
$K_{2\rho}^{-1}= K_{\a_1}^{-5} K_{\a_2}^{-8} K_{\a_3}^{-9}$.  If
we assume as in [Jn] that a short root satisfies $(\a,\a) = 2$, then
setting $(\e_i,\e_i) = 2$ for $i = 1, 2, 3$, we have that the
matrix of the bilinear form relative to the basis 
$\a_1 = \e_1 -\e_2$, $\a_2 = \e_2 - \e_3$, $\a_3 = \e_3$ of simple
roots is  

$$\left ( \matrix 4 & -2 & 0 \\ -2 & 4 & -2 \\ 0 & -2 & 2 \\ \endmatrix \right).$$

\n  From this we can compute that $(-2\rho, \e_i) = -10,-6,-2$ for $i = 1,2,3$,
respectively.    Since we know that for any weight vector $x_\mu$ that
$K_\a x_\mu = q^{(\mu,\a)}x_\mu$, we obtain 
that relative to the basis

$$\{x_{\e_1},x_{\e_2},x_{\e_3},x_0, x_{-\e_3},x_{-\e_2},x_{-\e_1}\},\tag 8.9$$

\n for the natural representation,

$$\phi_{\text{B}_3} =
\text{diag}\{q^{-10},q^{-6},q^{-2},1,q^2,q^6,q^{10}\}.\tag 8.10$$  

\n The spin representation of  $U_q($B$_3)$ has a basis 

$$\{x_{\xi_1},x_{\xi_2},x_{\xi_3},x_{\xi_4},x_{-\xi_4},x_{-\xi_3},x_{-\xi_2},x_{-\xi_1}\},$$

\n where $\xi_1 = (1/2)(\e_1+\e_2+\e_3)$, $\xi_2 =
(1/2)(\e_1+\e_2-\e_3)$, $\xi_3 = (1/2)(\e_1-\e_2+\e_3)$, and $\xi_4
= (1/2)(-\e_1+\e_2+\e_3)$.   The corresponding matrix of $\phi_{\text{B}_3}^+$
is given by

$$\phi_{\text{B}_3}^+  = \text{diag}\{q^{-9},q^{-7},q^{-3},q^{-1},
q,q^3,q^7,q^{9}\}.$$ 
\m
Now suppose that $f: U_q(\g_1) \rightarrow U_q(\g_2)$ is a Hopf algebra homomorphism
between two quantum groups.
Consider representations $(X,\pi)$ and $(X,\pi')$ for $U_q(\g_2)$ as above.
Then $(X,\pi \circ f)$, $(X,\pi' \circ f)$ are representations for
$U_q(\g_1)$.  Since $\pi' \circ f = \pi \circ S^2 \circ f =
\pi \circ f \circ S^2$, it follows that $\pi' \circ f = (\pi \circ f)'$.  
Moreover, the map

$$\phi_{\g_2}: (X,\pi) \rightarrow (X,\pi') $$

\n satisfies 

$$\phi_{\g_2}((\pi \circ f)(a)x) = \pi'(f(a)) \phi_{\g_2}(x) =
(\pi \circ f)'(a) \phi_{\g_2}(x)$$

\n for all $a \in U_q(\g_1)$ and $x \in X$ so it is an isomorphism
as $U_q(\g_1)$-modules:  $\phi_{\g_2}: (X, \pi \circ f) \rightarrow 
(X,(\pi \circ f)')$. 
\b
{\bf Proof of Proposition 8.4.}  Suppose $f: U_q($G$_2) \rightarrow
U_q($B$_3)$ is a Hopf algebra homomorphism.  Consider the natural
representation $Y = (X,\pi)$ of $U_q($B$_3)$, and the corresponding
representation

$$\phi_{\text{B}_3}: (X, \pi \circ f) \rightarrow (X, (\pi \circ f)')
\tag 8.11$$

\n of $U_q($G$_2)$-modules. Either these are the 7-dimensional irreducible
$U_q($G$_2)$-module or they are the sum of 7 trivial 1-dimensional modules. 
In the first case $\phi_{\text{B}_3}$ is a multiple of
$\phi_{\text{G}_2}$, so the eigenvalues
of $\phi_{\text{B}_3}$ are $\{\l q^{\pm 10},\l q^{\pm 8},\l q^{\pm 2},\l\}$.
This is impossible since $q$ is not a root of unity, so 
it must be that $\pi \circ f$ is trivial. 
 
Analogously, consider the spin representation $Y^+ = (X^+,\pi^+)$ of 
$U_q($B$_3)$.  Then 
as a $U_q$(G$_2)$-module, such a representation
should be the sum of representations of
dimension 1 or 7.  Reasoning as in (8.11), we see that
if a 7-dimensional $U_q$(G$_2)$-module occurs, then
$\phi_{\text{B}_3}^+$ restricted to the 7-dimensional submodule
should be proportional to $\phi_{\text{G}_2}$.  This would mean
that

$$\{\l q^{\pm 10},\l q^{\pm 8},\l q^{\pm 2},\l\} \subseteq 
\text{diag}\{q^{\pm 9},q^{\pm 7},q^{\pm 3},q^{\pm 1}\}$$

\n for some scalar $\l$, which is impossible.  Thus $\pi^+ \circ f$
is trivial.  
 
We have seen that the natural and spin representations for 
$U_q($B$_3)$ are trivial
modules for $f(U_q($G$_2))$.  Thus, $f(a) - \e(f(a))1$ is
in the annihilator of those representations. 
Since any finite-dimensional irreducible $U_q($B$_3)$-representation
of type I can be obtained from the tensor product of those, and since
$f$ is a Hopf homomorphism, it follows that $f(U_q($G$_2))$ acts
trivially on any finite-dimensional irreducible $U_q($B$_3)$-representation
of type I.   However,  $\cap\,\text{ann}(V) = (0)$ for all such representations
(see [Jn, Prop. 5.11]).  Consequently, $f(a) = \e(f(a))1$ for
all $a \in U_q($G$_2)$.   
\m
The other two cases in Proposition 8.4 are similar. For
a Hopf homomorphism $f: U_q(\text{G}_2) \rightarrow
U_q(\text{D}_4)$, it boils down to the fact that any finite-dimensional 
irreducible $U_q(\text{D}_4)$-representation occurs in tensor
products of copies of the natural or spin representations and 

$$\{\l q^{\pm 10}, \l q^{\pm 8}, \l q^{\pm 2}, \l\} \not \subseteq
\{q^{\pm 6},q^{\pm 4},q^{\pm 2},1\}.$$

\n In the final case,  $f: U_q(\text{B}_3) \rightarrow
U_q(\text{D}_4)$, the argument amounts to the impossibility of

$$\{\l q^{\pm 10}, \l q^{\pm 6}, \l q^{\pm 2}, \l\} \subseteq
\{q^{\pm 6},q^{\pm 4},q^{\pm 2},1 \}, \quad \text {and}$$

$$\{\l q^{\pm 9}, \l q^{\pm 7}, \l q^{\pm 3}, \l q^{\pm 1}\} \subseteq
\{q^{\pm 6},q^{\pm 4},q^{\pm 2},1 \}. \qed $$
\b
\head {\S 9. Quantum Octonions and Quantum Clifford Algebras} \endhead
\b
In this section we investigate the representation of  the quantum Clifford
algebra $C_q(8)$ on  $\Bbb O_q \oplus \Bbb O_q$.  Our main
reference for the quantum Clifford algebra is [DF], but we need to
make several adjustments to make the results of [DF] conform to
our work here.  It is convenient to begin with the general
setting of the quantum group of D$_n$ and then to specialize to D$_4$.  
\b
\subhead {Explicit form of $\check R_{V,V}$ for $U_q($D$_n)$} \endsubhead
\b
The natural representation $V$ of the quantum group $U_q($D$_n)$ has a basis

$$\{x_1, x_2, \dots, x_{2n}\} = \{x_{\e_1}, \dots, x_{\e_n}, x_{-\e_n},
\dots, x_{-\e_1}\}, \tag 9.1$$

\n where $E_\a, F_\a, K_\a$ act as in (2.2).  A straightforward computation
shows that relative to the basis $\{x_1, x_2, \dots, x_{2n}\}$, the
matrix of $\check R = \check R_{V,V}$ is given by

$$\aligned \check R & = q^{-1} \sum_{i = 1}^{2n} E_{i,i} \ot E_{i,i} 
+ \sum_{i \neq j,j'} E_{j,i} \ot E_{i,j} + (q^{-1} - q)
\sum_{j > i} E_{j,j} \ot E_{i,i}  \\
& \quad + q \sum_{i = 1}^{2n} E_{i',i} \ot E_{i,i'} + (q - q^{-1})
\sum_{i < j} (-q)^{\rho_{j'} - \rho_{i'}} E_{i',j} \ot E_{i,j'},
\\ \endaligned \tag 9.2$$

\n  where $i' = 2n+1 - i$ and $\rho = (\rho_1, \dots, \rho_{2n})
= (n-1,n-2,\dots,1,0,0,-1, \dots, 2-n, 1-n)$.  The projection 
$(\,,\,): V \ot V \rightarrow K$
is given by

$$(x_i, x_j) = \delta_{i,j'} (-1)^{n-1} (-q)^{\rho_{i'}}.$$

\n If we modify the basis slightly by defining

$$\aligned e_i =&  (-1)^{n-1+\rho_{i'}} x_i \\
e_{i'} =& x_{i'} \\
\endaligned \tag 9.3$$

\n for $i = 1, \dots, n$, then we obtain the usual expression for $\check R$:

$$\aligned \check R & = q^{-1} \sum_{i = 1}^{2n} E_{i,i} \ot E_{i,i} 
+ \sum_{i \neq j,j'} E_{j,i} \ot E_{i,j} + (q^{-1} - q)
\sum_{j > i} E_{j,j} \ot E_{i,i}  \\
& \quad + q \sum_{i = 1}^{2n} E_{i',i} \ot E_{i,i'} + (q - q^{-1})
\sum_{i < j} q^{(\rho_{j'} - \rho_{i'})} E_{i',j} \ot E_{i,j'},
\\ \endaligned \tag 9.4$$

\n and the bilinear form is given by

$$(e_i, e_j) = \delta_{i,j'}q^{\rho_{i'}}.$$

\m
\subhead {The quantum Clifford algebra} \endsubhead
\b
In [DF], Ding and Frenkel define the quantum Clifford algebra $C_q(2n)$
as the quotient of the tensor algebra $T(V)$ of $V$ by the ideal generated
by $\{v \ot w + q \check R'_{2,1}(v \ot w) - (w,v)1\mid v,w \in V\}$ where $\check R'$
is $\check R$ with $q$ and $q^{-1}$ interchanged, and
$\check R'_{2,1} = \sigma \circ \check R' \circ \sigma$ where
$\sigma(v \ot w) = (w \ot v)$.   This quantum Clifford algebra can be thought
of as a unital associative algebra generated by $\psi_1, \dots, \psi_n,
\psi_1^*, \dots, \psi_n^*$ with relations,

$$\aligned & \psi_i \psi_j = -q^{-1}\psi_j \psi_i \quad 
\quad \psi_i^* \psi_j^* = -q\psi_j^* \psi_i^* \quad (i > j) \\
& \psi_i \psi_i =  0 = \psi_i^* \psi_i^* \quad \quad
\psi_i \psi_j^* =  -q^{-1}\psi_j^* \psi_i \quad (i \neq j) \\
&\psi_i \psi_i^* + \psi_i^* \psi_i = (q^{-2} -1) \sum_{i < j} \psi_j \psi_j^* +1.\\
\endaligned \tag 9.5$$    

If we make the change of notation
$Y_i = \psi^*_{n-i}, Y_i^* = \psi_{n-i}$, $ i = 1, \dots, n$, then the
relations in (9.5) become:

$$\aligned & Y_i Y_j = -q^{-1}Y_j Y_i \quad 
\quad Y_i^* Y_j^* = -q Y_j^* Y_i^* \quad (i > j) \\
& Y_i Y_i =  0 = Y_i^* Y_i^*  \quad \quad
Y_i Y_j^* =  -q Y_j^* Y_i \quad (i \neq j) \\
&Y_i Y_i^* + Y_i^* Y_i = (q^{-2} -1) \sum_{i > j} Y_j Y_j^* +1.\\
\endaligned \tag 9.6$$

In order to construct a representation of $C_q(8)$ on $\Bbb O_q \oplus \Bbb O_q$ 
we will need the following proposition: 
\b
\proclaim{Proposition 9.7}  Let $I$ be the ideal of the tensor 
algebra $T(V)$ generated by
$\{v \ot w + q^{-1} \check R_{2,1}(v \ot w) - (w,v)1\mid v,w \in V\}$
where $\check R_{2,1} = \sigma \circ \check R \circ \sigma$.  Then $C_q(2n)
\cong T(V)/I$.    \endproclaim
\b
{\bf Proof.} Consider the following basis of $V$:

$$\{y_i = q^{\rho_i}e_i,\ \ y_i^* = e_{i'}, \mid i = 1, \dots, n\}\tag 9.8$$

\n which satisfies

$$(y_i,y_j^*) = \delta_{i,j} \qquad \qquad  (y_j^*,y_i) = \delta_{i,j}q^{2\rho_i}.
\tag 9.9$$  

\n Let $g: V \ot V \rightarrow I$ be defined by  

$$ g(v \ot w) = v \ot w + q^{-1} \check R_{2,1}(v \ot w) - (w,v)1. \tag 9.10$$

\n From the explicit form of $\check R$ on the basis $\{e_i\}$, one can check that
the following hold:

$$\aligned & g(y_i \ot y_j) = y_i \ot y_j +q^{-1}y_j \ot y_i \quad  (i > j) \\
& g(y_j^*\ot y_i^*) = y_j^* \ot y_i^* + q^{-1}y_i^* \ot y_j^* \quad (i > j) \\
& g(y_i\ot y_i) = (1+q^{-2})y_i \ot y_i \quad \quad g(y_i^* \ot y_i^*)
= (1+q^{-2})y_i^* \ot y_i^* \\
& g(y_j^* \ot y_i) =  y_j^* \ot y_i + q^{-1}y_i \ot y_j^*  \quad (i \neq j) \\
&g(y_i^* \ot y_i) =  y_i^* \ot y_i + y_i \ot y_i^* + (1-q^{-2}) \sum_{i > j} y_j\ot y_j^*
-1, 
\endaligned \tag 9.11$$  

\n These give the relations in (9.6) in $T(V)/I$.   However, we also have the
relations:

$$\aligned & g(y_j \ot y_i) = q^{-1}y_i \ot y_j +q^{-2}y_j \ot y_i \quad (i > j) \\ 
& g(y_i^*\ot y_j^*) = q^{-1}y_j^* \ot y_i^* + q^{-2}y_i^* \ot y_j^* \quad (i > j) \\ 
& g(y_i \ot y_j^*) =  q^{-1}y_j^* \ot y_i + q^{-2}y_i \ot y_j^*  \quad (i \neq j) \\
& g(y_i \ot y_i^*)
= y_i^* \ot y_i +  q^{-2} y_i \ot y_i^* \\
&\hskip 1 truein  + (1-q^{-2})q^{2\rho_i} 
\Big(\sum_{j = 1}^n y_j\ot y_j^* + \sum_{j > i}q^{-2\rho_j} y_j^* \ot y_j\Big)
-q^{2\rho_i}. 
\endaligned \tag 9.12$$  

\n The first three identities don't say anything new.  So to conclude the
proposition, it suffices to verify that the relation 
$g(y_i \ot y_i^*) = 0$ holds in $C_q(2n)$, that is:

\b
\proclaim{Lemma 9.13}  In $C_q(2n)$ we have 

$$ q^{-2} Y_iY_i^* + Y_i^* Y_i + (1-q^{-2})q^{2\rho_i}
\Big(\sum_{j = 1}^n Y_jY_j^* + \sum_{j > i}q^{-2\rho_j} Y_j^*Y_j\Big) 
= q^{2\rho_i}.$$
\endproclaim 
\b
{\bf Proof.}  As a shorthand notation we set $a_i = Y_iY_i^*$, 
$a_i' = Y_i^* Y_i$, and $b_i = a_i + a_i'$ for $i = 1, \dots, n$.   From the relations in
$C_q(2n)$ we know that $b_i = (q^{-2}-1)\sum_{j < i}a_j + 1$ so
that $b_{i+1} = (q^{-2}-1)a_i + b_i$ or 

$$a_{i+1}+a_{i+1}'
= q^{-2}a_i + a_i'. \tag 9.14 $$   

\n Letting  $c_i = q^{-2}a_i + a_i'$  we have 

$$c_i -q^2 c_{i+1} = q^{-2}a_i + a_i' -a_{i+1} -q^2 a_{i+1}' \eqa
 (1-q^2) a_{i+1}'$$

\n where (1) follows from (9.14).   Therefore,

$$\aligned c_i & = (1-q^2)a_{i+1}' + q^2(1-q^2)a_{i+2}' + \cdots + q^{2(n-i-1)}(1-q^2)a_n'
+ q^{2(n-i)}c_n \\
& = (q^{-2}-1)q^{2\rho_i}\Big(\sum_{j>i}q^{-2\rho_j}a_j'\Big) + q^{2\rho_i}c_n.\\
\endaligned \tag 9.15$$

\n On the other hand,

$$\aligned c_n & = q^{-2}a_n + a_n' = (q^{-2}-1)a_n + a_n + a_n' = (q^{-2}-1)a_n + b_n \\
& = (q^{-2}-1)a_n + (q^{-2}-1)\Big(\sum_{j < n}a_j\Big) + 1 = (q^{-2}-1)
\Big(\sum_{j \leq n} a_j\Big) +
1.
\\
\endaligned 
\tag 9.16$$

\n The left-hand side of the identity in the statement of the lemma is

$$\aligned c_i + (1-q^{-2})q^{-2\rho_i}\Big( \sum_{j = 1}^n a_j\Big) & + 
(1-q^{-2})q^{2\rho_i}\Big(\sum_{j > i}q^{-2\rho_j}a_j'\Big) \\
& \eqa c_i + q^{2\rho_i}(1-c_n) + (1-q^{-2})q^{2\rho_i}\Big(\sum_{j > i}q^{-2\rho_j}a_j'\Big)
\\ & \eqb q^{2\rho_i} \quad \text{(the right side)},\\ \endaligned $$

\n where (1) is (9.16) and (2) is (9.15).  \qed

\b
\subhead {The $\uqd$ case} \endsubhead
\b
Let us return now to our quantum octonions.   Consider the
two ``intertwining operators'' $\Phi^0$
and $\Phi^1$ which are the $\uqd$-module homomorphisms:

$$
\aligned & \Phi^0: V_\eta \ot V \jimap V_{\eta \circ \zeta} \ot V \jpmap V_\zeta \\
&\Phi^1: V_\zeta \ot V \jimap V_{\zeta \circ \eta} \ot V \pmap V_\eta, \\ 
\endaligned $$  

\n We scale these mappings by setting $\Psi^0 = q^{-\frac32} \Phi^0$ and $\Psi^1 = 
q^{-\frac32} \Phi^1$.
Let $\frak S_q = ($id$ + q^{-1}\check R)$, where $\check R =\check R_{V,V}$, which
we know maps $V\ot V$ onto ${\Cal S}^2(V)$. 
\b
\proclaim{Lemma 9.17} (Compare [DF, Prop. 3.1.2])

$$
\aligned
& \Psi^1 \circ (\Psi^0 \ot \text{id})\circ ( \text{id} \ot \frak S_q) = 
 \text{id} \ot (\,,\,) \\
&  \Psi^0 \circ (\Psi^1 \ot \text{id})\circ ( \text{id} \ot \frak S_q) = 
 \text{id} \ot (\,,\,). \\ \endaligned
$$  \endproclaim
\b
{\bf Proof.}  It is enough to check the values on $x_1 \ot x_1 \ot  x_{-1}$:

$$
\aligned \Psi^1 \circ (\Psi^0 \ot \text{id})\circ ( \text{id} \ot \frak S_q)(
x_1 \ot x_1 \ot  x_{-1}) & = \Psi^1 \circ (\Psi^0 \ot \text{id})(x_1
\ot (x_1 \ot x_{-1} + x_{-1} \ot x_1)) \\
& = q^{-\frac 32}\Psi^1(q^{-\frac 32} x_4 \ot x_1) = q^{-3} x_1. \\
\endaligned $$ 

\n But $(x_1,x_{-1}) = (q^3 + q^{-3})(x_1 | x_{-1}) = q^{-3}$, so both
maps agree on $x_1 \ot x_1 \ot  x_{-1}$.  The computation for the
other map is left as an exercise.  \qed
\b
Now for any $v \in V$ define

$$\aligned \Psi(v): & V_\zeta \oplus V_\eta \rightarrow V_\zeta \oplus V_\eta \\
 & \Psi(v)(x+y) = \Psi^0(y \ot v) + \Psi^1(x \ot v). \\
\endaligned \tag 9.18$$

\n The relations in the previous lemma say that

$$\Psi \circ (\Psi \ot \text{id}) \circ (\text{id} \ot \frak S_q)
= \text{id}\ot (\,,\,),$$

\n so for any basis elements $e_i, e_j \in V$ and $x \in V_\zeta \oplus
V_\eta$, we have

$$
\aligned
(e_j,e_i)x & = \Psi \circ (\Psi \ot \text{id}) \circ (\text{id} \ot \frak S_q)(x \ot e_j
\ot e_i) \\
& = \Psi \circ (\Psi \ot \text{id}) (x \ot (e_j \ot e_i + q^{-1} \check R(e_j \ot e_i))) \\
& = \Psi(e_i) \Psi(e_j)(x) + q^{-1}\sum_{k,\ell}
(\check R \circ \sigma)^{k,\ell}_{i,j} \Psi(e_\ell)\Psi(e_k)(x), \\
\endaligned$$

\n where $(\check R \circ \sigma)^{k,\ell}_{i,j}$ is the coordinate matrix of
$\check R \circ \sigma$ in the basis $\{e_i \ot e_j \}$.   So
$\Psi(e_i)$ and $\Psi(e_j)$ satisfy

$$(e_j,e_i)\text{id} = \Psi(e_i)\Psi(e_j) + q^{-1}
\sum_{k,\ell}(\check R \circ \sigma)^{k,\ell}_{i,j}\Psi(e_\ell)\Psi(e_k).$$

\n Notice now that

$$\aligned
e_i \ot e_j + q^{-1}\check R_{2,1}(e_i \ot e_j) & = e_i \ot e_j+
q^{-1} \sigma \circ \check R \circ \sigma (e_i \ot e_j) \\
& = e_i \ot e_j + q^{-1} 
\sum_{k,\ell}(\check R \circ \sigma)^{k,\ell}_{i,j} e_\ell \ot e_k,\\
\endaligned
$$

\n so the operators $\Psi(v)$ verify the relations of the quantum Clifford
algebra.   This shows that the map

$$\aligned \Psi: V &  \rightarrow \End(V_\zeta \oplus V_\eta) \\
v & \mapsto \Psi(v): x + y \rightarrow q^{-3}\Big(\jmath(y \cdot \jmath(v)) + \jmath(x)
\cdot v\Big)
\\ \endaligned \tag 9.19$$

\n extends to a representation $\Psi: C_q(8) \rightarrow \End(V_\zeta \oplus V_\eta)
= \End(\Bbb O_q \oplus \Bbb O_q)$.  
\b
In [DF], Ding and Frenkel show that $C_q(2n)$ is isomorphic to the classical Clifford
algebra, so $C_q(2n)$ is a central simple associative algebra of dimension $2^8$.  
The same is true of $\End(\Bbb O_q \oplus \Bbb O_q)$.  As a result we have
\b
\proclaim {Proposition 9.20} 
The map $\Psi: C_q(8) \rightarrow \End(V_\zeta \oplus V_\eta)
= \End(\Bbb O_q \oplus \Bbb O_q)$ in (9.19) is an algebra isomorphism.  
\endproclaim 

\b
It is worth noting that at $q = 1$ the map in Proposition 9.20 reduces to 

$$\aligned
\Psi:& \Bbb O \rightarrow \End(\Bbb O \oplus \Bbb O)  \\
& v \mapsto \Psi(v): x + y \mapsto  
v \cdot \overline y  + \overline x \cdot v \\
\endaligned \tag 9.21 $$

\n so 

$$\Psi(v)\Psi(v) : x + y \mapsto \Psi(v)(v \cdot \overline y +\overline x \cdot v)
= v \cdot (\overline v \cdot x) + (y \cdot \overline v)\cdot v   = (v|v)_{{}_1}(x + y),$$

\n (see Props. 2.21, 4.3).  This shows that
$\Psi$ extends to the standard Clifford algebra - that is, to the  
quotient $T(\Bbb O)/ < v \ot v - (v|v)_{{}_1}1>$  of the tensor algebra,
in this case as in [KPS]. 
\vskip .5 truein  

\def\m{

\enddocument
\end